\documentclass[10pt]{article}
\usepackage{hyperref}
\usepackage{amssymb}
\usepackage{amsmath}
\input{liemacs10.sty} 
\newcommand{\BGL}{\mathop{{\rm BGL}}\nolimits}

\renewcommand{\phi}{\varphi}

\newcommand{\Stand}{\mathop{{\rm Stand}}\nolimits}

\renewcommand\mlabel{\label}

\begin{document}

\title{Finite dimensional semigroups\\ of unitary endomorphisms of standard subspaces}
\author{Karl-H. Neeb\begin{footnote}{
Department  Mathematik, FAU Erlangen-N\"urnberg, Cauerstrasse 11, 
91058-Erlangen, Germany; neeb@math.fau.de}
\end{footnote}
\begin{footnote}{Supported by DFG-grant NE 413/9-1.}
\end{footnote}}

\maketitle

\abstract{Let $V$ be a standard subspace in the complex 
Hilbert space $\cH$ and $G$ be a finite dimensional Lie  
group of unitary and antiunitary operators on $\cH$ 
containing the modular group $(\Delta_V^{it})_{t \in \R}$ 
of $V$ and the corresponding modular conjugation~$J_V$. 
We study the semigroup 
\[ S_V = \{ g\in G \cap \U(\cH)\: gv \subeq V\} \] 
and determine its Lie wedge $\L(S_V) = \{ x \in \g \: 
\exp(\R_+ x) \subeq S_V\}$, i.e., the generators of its one-parameter 
subsemigroups in the Lie algebra $\g$ of~$G$. 
The semigroup $S_V$ is analyzed in terms 
of antiunitary representations and their analytic extension 
to semigroups of the form $G \exp(iC)$, where $C \subeq \g$ is an 
$\Ad(G)$-invariant closed convex cone. 

Our main results assert that 
the Lie wedge $\L(S_V)$ spans a $3$-graded Lie subalgebra 
in which it can be described explicitly in terms 
of the involution $\tau$ of $\g$ induced by $J_V$, 
the generator $h \in \g^\tau$ of the modular group, 
and the positive cone of the corresponding representation. We also derive 
some global information on the semigroup $S_V$ itself.}

\tableofcontents 

\vspace{1cm}

\section{Introduction}
\mlabel{sec:1}

Let $\cH$ be a complex Hilbert space and 
$\cM \subeq B(\cH)$ be a von Neumann algebra. 
Further, let $\Omega \in \cH$ be a unit vector which is 
{\it cyclic} for $\cM$ ($\cM\Omega$ is dense in $\cH$) and 
{\it separating} (the map $\cM \to \cH,M \mapsto M\Omega$ is injective). 
By the Tomita--Takesaki Theorem (\cite[Thm.~2.5.14]{BR87}), 
the closed real subspace $V := V_\cM := \oline{\{ M\Omega \: M = M^* \in \cM\}}$ 
is {\it standard}, i.e., 
\begin{equation}
  \label{eq:stansub}
  V \cap i V = \{0\}\quad \mbox{ and } \quad \cH = \oline{V + i V}
\end{equation} 
(cf.\ \cite{Lo08} for the basic theory of standard subspaces). 
To the standard subspace $V$, we can associate a {\it pair of modular 
objects} $(\Delta_V, J_V)$, i.e., $\Delta_V > 0$ is a positive 
selfadjoint operator, $J_V$ is a {\it conjugation} (an antiunitary 
involution), and these two operators satisfy the modular relation 
$J_V \Delta_V J_V = \Delta_V^{-1}$. The pair $(\Delta_V, J_V)$ is obtained
by the polar decomposition $\sigma_V = J_V \Delta_V^{1/2}$ of the closed operator 
\[ \sigma_V \: \cD(\sigma_V) := V + i V \to \cH, \quad 
x + i y \mapsto x- iy \] 
with $V = \Fix(\sigma_V)$. The main assertion of the 
Tomita--Takesaki Theorem is  that 
\[ J_V \cM J_V = \cM' \quad \mbox{ and } \quad 
\Delta_V^{it} \cM \Delta_V^{-it}  = \cM \quad \mbox{ for }  \quad t \in \R.\] 
So we obtain a one-parameter group of automorphisms of $\cM$ 
(the modular group) 
and a symmetry between $\cM$ and its commutant $\cM'$, implemented by $J_V$. 

Motivated by the Haag--Kastler theory 
of local observables in Quantum Field Theory (QFT) 
(\cite{Ha96}, \cite{BS93}, \cite{BDFS00}), 
we are interested in finite dimensional Lie groups 
$G \subeq \U(\cH)$ of unitary operators fixing $\Omega$, containing 
the corresponding modular group $(\Delta_V^{it})_{t \in \R}$ 
and invariant under conjugation with the modular conjugation~$J_V$. 
In this context, we would like to understand the subsemigroup 
\[ S_\cM := \{ g \in G \: g \cM g^{-1} \subeq \cM\} \] 
of those elements of $G$ acting by endomorphisms on $\cM$ 
(\cite{BLS11, LL15, Le15}). 
As $g\Omega = \Omega$ for $g \in G$, we have
$V_{g\cM g^{-1}} = g V_\cM$, so that 
$g\cM g^{-1} \subeq \cM$ implies $gV_\cM \subeq V_\cM$. 
For $V = V_\cM$, we therefore have 
\begin{equation}
  \label{eq:SM=SV}
S_\cM \subeq S_V :=  \{ g\in G \: gV \subeq V\}.
\end{equation}
It follows in particular that, if $S_\cM$ has interior points, 
then so does the semigroup $S_V$. 
In the present paper we determine its {\it Lie wedge} 
\begin{footnote}{In the theory of 
Lie semigroups (\cite{HHL89, HN93}) Lie wedges are the semigroup analogs 
of the Lie algebras of  closed subgroups. A {\it Lie wedge} is a closed convex 
cone $W$ in a Lie algebra $\g$ such that 
$e^{\ad x} W = W$ for $x \in W \cap -W$. In particular, 
linear subspaces are Lie wedges if and only if they are Lie subalgebras.}
\end{footnote}
\[ \L(S_V) = \{ x \in \g \: 
\exp(\R_+ x) \subeq S_V\}, \] 
i.e., the set of generators of its one-parameter 
subsemigroups in the Lie algebra $\g$ of $G$ (\cite{HHL89, HN93}). 

The current interest in standard subspaces arose 
in the 1990s from the work of Borchers and Wiesbrock (\cite{Bo92, Wi93}). 
This in turn led to the concept of modular localization 
in Quantum Field Theory introduced by  
Brunetti, Guido and Longo in \cite{BGL02, BGL94, BGL93}. 
We refer to Subsection~\ref{subsec:5.0} for more on the relation 
to von Neumann algebras. 

Compared to the rather inaccessible object $S_\cM$, 
the semigroup  $S_V$ can be analyzed in terms of {\it antiunitary representations} 
of {\it graded Lie groups}: A graded Lie group is a pair 
$(G,\eps_G)$, where $\eps_G \: G \to \{\pm 1\}$ is a homomorphism 
and we write $G_{\pm} = \eps_G^{-1}(\pm 1)$, so that $G_+ \trile G$ is a 
normal subgroup of index $2$ and $G_- = G \setminus G_+$. 
An important example is the group 
$\AU(\cH)$ of unitary or antiunitary operators on a complex 
Hilbert space with $\AU(\cH)_+ =\U(\cH)$. A morphism of graded groups 
$U \: G \to \AU(\cH)$ is called an antiunitary representation. 
Then $U(G_+) \subeq \U(\cH)$ and $U(G_-)$ consist of antiunitary operators. 

We write $\Stand(\cH)$ for the set of standard subspaces of~$\cH$. 
We have already seen that every standard subspace 
$V$ determines a pair $(\Delta_V, J_V)$ of modular objects 
and that $V$ can be recovered from this pair by $V = \Fix(J_V \Delta_V^{1/2})$. 
This observation can be used to obtain a representation theoretic 
parametrization of $\Stand(\cH)$: 
each standard subspace $V$ specifies a homomorphism 
\begin{equation}
  \label{eq:uv-rep}
 U^V \: \R^\times \to \AU(\cH)\quad \mbox{ by } \quad 
U^V(e^t) := \Delta_V^{-it/2\pi}, \quad 
U^V(-1) := J_V.
\end{equation}
We thus obtain a bijection between $\Stand(\cH)$ and 
antiunitary representations of the graded Lie group $\R^\times$ 
with $\eps(r) = \sgn(r)$ (\cite{NO17}). 
For a given antiunitary representation $(U,\cH)$ of a graded Lie group 
$(G,\eps_G)$, we thus obtain a natural map, 
the {\it Brunetti--Guido--Longo map} 
\begin{equation}
  \label{eq:bgl}
\BGL \: \Hom_{\rm gr}(\R^\times,G) \to \Stand(\cH), \quad 
\gamma \mapsto V_\gamma \quad \mbox{ with } \quad   
U^{V_\gamma} = U \circ \gamma
\end{equation}
(\cite{BGL02}, \cite{NO17}). 
Note that $\gamma \in \Hom_{\rm gr}(\R^\times,G)$ is completely determined by 
\[ h := \gamma'(1) \in \g \quad \mbox{ and } \quad \sigma := \gamma(-1).\] 
As $\sigma^2 = e$, it defines an involution $\tau_G(g) := \sigma g \sigma$ 
on~$G$, an involution $\tau = \Ad(\sigma)$ on~$\g$ with $\tau(h) = h$, and 
$G \cong G_+ \rtimes \{\id_G, \tau_G\}$. 

We thus arrive at the problem to determine 
for an injective antiunitary representation $(U,\cH)$ of a graded Lie group 
$(G, \eps_G)$ and a standard subspace $V = V_\gamma \subeq \cH$ obtained by the 
BGL construction from a pair $(\tau,h)$, consisting 
of an involutive automorphism $\tau$ of $\g$ and an element 
$h \in \g$ with $\tau(h) = h$, the semigroup 
\[ S_V := \{ g \in G_+ \: U(g) V \subeq V \}.\] 
A crucial piece of information on $S_V$ 
is contained in its Lie wedge $\L(S_V)$. 
To formulate our main results, for $\lambda \in \R$ and 
an $\ad h$-invariant subspace $F \subeq \g$, 
write $F_\lambda(h) := \ker(\ad h - \lambda \id_\g) \cap  F$ for the 
corresponding eigenspace. We also put 
$\fh := \ker(\tau - \id)$ and $\fq := \ker(\tau + \id)$ 
and write $C_U := \{ x \in \g \: - i \partial U(x) \geq 0\}$ for the 
{\it positive cone of $U$}.
The Structure Theorem (Theorem~\ref{thm:main2}) asserts that 
\begin{equation}
  \label{eq:(b)}
\L(S_V) = C_- \oplus \fh_0(h) \oplus C_+ 
\end{equation}
for the two pointed closed convex cones 
\[ C_\pm := \L(S_V) \cap \fq_{\pm 1}(h) = \pm C_U \cap \fq_{\pm 1}(h). \] 
Further, $\L(S_V)$ spans 
a $3$-graded Lie subalgebra $\g_{\rm red}$, and the 
cones $C_\pm$ are abelian subsets of~$\g$. 

So we obtain an explicit description 
of the Lie wedge $\L(S_V)$ in terms of 
the positive cone $C_U$ of the representation $(U,\cH)$, 
the involution $\tau$ of $\g$ induced by $J_V$, 
and the generator $h \in \g^\tau$ of the modular group. 
It shows in particular that the most interesting situations are those 
where $\g$ is $3$-graded by $\ad h$, i.e., 
$\g = \g_{-1}(h) \oplus \g_0(h) \oplus \g_1(h)$, and 
$\tau = e^{\pi i \ad h}$. In this context, the representation $U$ should be such that 
the cones $C_U \cap \fg_{\pm 1}(h)$ generate $\g_{\pm 1}(h)$. 
We refer to Subsection~\ref{subsec:5.1a} for more comments on 
classification problems.

One of our key tools is a characterization of the operators contained 
in the algebra 
$\cA_V := \{ A \in B(\cH) \: A V \subeq V\}$ of {\it $V$-real operators} 
in terms of the orbit maps 
\[ \alpha^A(t) := \alpha_t(A) := \Delta_V^{-it/2\pi} A \Delta_V^{it/2\pi}.\]
The Araki--Zsid\'o Theorem (\cite{AZ05}) asserts that, for $A \in B(\cH)$,  
$A \in \cA_V$ is equivalent to the existence of an analytic continuation 
of $\alpha^A$ from $\R$ to the closure of the strip 
\[ \cS_\pi 
= \{ z \in \C \: 0 <  \Im z < \pi\}\] satisfying 
$\alpha^A(\pi i) = J_V A J_V$. 
It follows in particular, that $\cA_V$ is invariant under the 
involution $A^\sharp := J_V A^* J_V$,  and that we 
obtain for every $z \in \oline{\cS_\pi}$ an injective 
representation 
\[ \alpha_z \: \cA_V \to B(\cH), A \mapsto \alpha^A(z) \quad 
\mbox{ with } \quad \|\alpha_z\| \leq 1.\]
For $z = \frac{\pi i}{2}$ we even obtain a 
$*$-representation $\alpha_{\frac{\pi i}{2}} \: \cA_V \to B(\cH^{J_V}), 
A \mapsto \hat A$ by operators commuting with~$J_V$. 

On the Lie group side, we mimic the Araki--Zsid\'o Theorem 
as follows. For a unitary representation
 $U \: G \to U(\cH)$ of a Lie group $G$, assumed with discrete 
kernel, we can extend $U$ to a representation of a  semigroup 
\[ S^U =  G \exp(iC_U), \] 
where $C_U$ is the positive cone of $U$, 
and the polar map $G \times C_U \to S^U, (g,x) \mapsto g \exp(ix)$ 
is a homeomorphism.\begin{footnote}
{Such semigroups are called 
{\it Olshanski semigroups}. They  first appear in 
Olshanski's paper \cite{Ol82} and an exposition of their 
theory can be found in \cite{Ne00}. The refinements needed 
for representations with non-discrete kernel have recently been worked 
out in \cite{Oeh18}. 
}\end{footnote}
Then $U(g \exp(ix)) = U(g) e^{i \partial U(x)}$ provides an extension 
of $U$ to $S^U$ (\cite[\S XI.2]{Ne00}). 
To bring modular conjugations into the picture, we also  consider 
an involution $\tau_G \in \Aut(G)$ (inducing an involution $\tau$ on $\g$),  
for which $U$ extends 
to an antiunitary representation of the graded Lie group 
$G \rtimes \{\id_G,\tau_G\}$. 
Then $J := U(\tau_G)$ is a conjugation satisfying 
$U(\tau_G(g)) = J U(g) J$ for $g \in G$. 
For $h \in \fh= \g^\tau$ we can now consider the 
standard subspace $V$ determined by $J_V = J$ 
and $\Delta_V = e^{2\pi i \partial U(h)}$, so that $\Delta_V^{-it/2\pi} = U(\exp th)$.  
Now the role  of the semigroup $\cA_V \cap \U(\cH)$ 
in the Araki--Zsid\'o Theorem is played by the subsemigroup 
$S^U_{\rm inv} \subeq G$, consisting of all elements 
$s \in G$ for which the orbit map $\beta^s(t) 
= \exp(th) s \exp(-th)$, defined on $\R$,  
extends analytically to a map from the closure of the strip 
$\cS_\pi$ to $S^U$, such that $\beta_{\pi i}(s) = \tau_G(s)$. 
Our second main result is the Inclusion Theorem 
(Theorem~\ref{thm:4.17}), asserting that $S^U_{\rm inv}\subeq S_V$. 
It is used to obtain one inclusion in the Structure 
Theorem mentioned above. It has a partial converse in the 
Germ Theorem (Theorem~\ref{thm:keylem3}) 
which shows that both subsemigroup have the same germ, 
i.e., that there exists an $e$-neighborhood $\cU\subeq G$ with 
$S^U_{\rm inv} \cap \cU= S_V \cap \cU$. 

The content of this paper is as follows. 
In Section~\ref{sec:2} we study for a standard subspace 
$V \subeq \cH$ the semigroup 
$S_V = \{ g \in \U(\cH) \: gV \subeq V\}$ of all 
unitary endomorphisms of~$V$. First we observe that 
$S_V$ is a group if and only if $\Delta_V$ is bounded, 
so that the situation is only 
interesting if $\Delta_V$ is unbounded (Lemma~\ref{lem:3.2}). 
We also state the Araki--Zsid\'o Theorem 
(a complete 
proof is provided in Appendix~\ref{app:a}) and develop 
its consequences. 

In Section~\ref{sec:3} we prepare the ground 
for our analysis of the subsemigroup $S^U_{\rm inv} \subeq G \subeq S^U$ 
which provides a Lie theoretic framework for verifying the Araki--Zsid\'o 
condition. The main result in Section~\ref{sec:3} is 
the Inclusion Theorem $S^U_{\rm inv}\subeq S_V$ (Theorem~\ref{thm:4.17}). 
Since both semigroups $S^U_{\rm inv}$ and $S_V$ are hard to 
describe globally, an important consequence of the 
Inclusion Theorem is the inclusion $\L(S^U_{\rm inv}) \subeq \L(S_V)$. 
To use this inclusion to prove the Structure Theorem, we 
derive an explicit description of the wedge 
$\L(S^U_{\rm inv})$ by interpreting it as a similar object 
$\L(S^U)_{\rm inv} = (\g + i C_U)_{\rm inv}$ in the abelian context. 

This motivates our independent discussion of the case where 
$G$ is a real Banach space $E$, endowed with an 
involution $\tau$ and an operator $h\in B(E)$, and 
$W \subeq E$ is a pointed closed convex cone invariant 
under $-\tau$ and the one-parameter group $e^{\R h}$ 
(Subsection~\ref{subsec:3.1}). 
In this simple situation the semigroup $(E + i W)_{\rm inv}$ 
can be determined very explicitly by elementary means and provides an 
important prototype for the more general non-abelian situation: 
\[ (E + i W)_{\rm inv} = (W \cap E^-_1(h)) \oplus E^+_0(h)
\oplus (-W \cap E^-_{-1}(h)), 
\quad \mbox{ where } \quad 
E^\pm = \ker(\tau \mp \1).\]
In Subsection~\ref{subsec:3.2} we then recall 
the basic facts on Olshanski semigroups $\Gamma_G(W) = G \exp(iW)$ 
for invariant cones $W \subeq \g$. They are non-abelian generalizations 
of the tubes $E + i W$. We prove the Inclusion 
Theorem in Subsection~\ref{subsec:3.3} and by applying it to 
the corresponding 
Lie wedges, we already obtain one inclusion of the Structure Theorem. 

The proof of the Structure Theorem 
(Theorem~\ref{thm:main2}) is completed in Subsection~\ref{subsec:4.1}, 
where we also prove the Germ Theorem. 
In Subsection~\ref{subsec:4.2} we describe the unit group 
$G_V$ of $S_V$, and in Subsection~\ref{subsec:4.3} 
we discuss some classes of examples. 
We conclude this paper with Section~\ref{sec:5} on perspectives and open 
problems. Some results that we did not find in the appropriate form 
in the literature are stated and proved in appendices. 

\subsection*{Notation} 

\begin{itemize}
\item For a Lie group $G$, we write $\g$ for its Lie algebra, 
$\Ad \: G \to \Aut(\g)$ for the adjoint action of $G$ on $\g$, induced by the 
conjugation action of $G$ on $G$, and $\ad x(y) = [x,y]$ for the adjoint 
action of $\g$ on itself. 
\item $(G,\eps_G)$ denotes a graded group, where 
$\eps_G \: G \to \{\pm 1\}$ is a homomorphism; 
$G_{\pm} = \eps_G^{-1}(\pm 1)$. An important example is the group 
$\AU(\cH)$ of unitary or antiunitary operators on a complex 
Hilbert space $\cH$ with $\AU(\cH)_+ =\U(\cH)$. A morphism of graded groups 
$U \: G \to \AU(\cH)$ is called an antiunitary representation. 
If $G$ is a topological group, then antiunitary representations are 
assumed to be continuous with respect to the strong operator topology on 
$\AU(\cH)$. 
\item For a graded homomorphism $\gamma \: \R^\times \to G$, we write 
$\sigma := \gamma(-1)$, $\tau = \Ad(\sigma)$, and $h := \gamma'(1) \in \g^\tau$. 
Then $\g = \fh \oplus \fq$ 
for the $\tau$-eigenspaces $\fh = \ker(\tau - \id_\g)$ and 
$\fq = \ker(\tau + \id_\g)$. 
We further write 
$\tau_G(g) := \sigma g \sigma$ for the corresponding involution 
on $G$. 
\item For a real standard subspace $V \subeq \cH$, we write 
$(\Delta_V, J_V)$ for the corresponding pair 
of {\it modular objects} with $V = \Fix(J_V \Delta_V^{1/2})$. 
\item Horizontal strips in the complex plane are denoted 
$\cS_{\alpha,\beta} := \{ z \in \C \: \alpha < \Im z < \beta\}$ and 
we also abbreviate $\cS_\beta := \cS_{0,\beta}$ for $\beta > 0$. 
\item For a unitary representation $U \: G \to \U(\cH)$ of a finite 
dimensional Lie group~$G$, 
we write $\cH^\infty$ for the dense subspace of {\it smooth vectors} $\xi$, for which 
the orbit maps $U^\xi \: G \to \cH, g \mapsto U_g \xi$ is smooth. 
We also have the dense subspace $\cH^\omega \subeq \cH^\infty$ of 
{\it analytic vectors} for which the orbit map $U^\xi$ is analytic. 
On $\cH^\infty$ 
we have a representation $\dd U$ of the complex Lie algebra $\g_\C$ 
given on $x \in \g$ by $\dd U(x)\xi = \derat0 U(\exp tx)\xi$. 
The infinitesimal generator of the unitary one-parameter group 
$(U(\exp tx))_{t \in \R}$ is denoted $\partial U(x)$. It coincides 
with the closure of the operator $\dd U(x)$. 
The closed convex $\Ad(G)$-invariant cone  
\[ C_U := \{ x \in \g \: -i \partial U(x) \geq 0 \} \] 
is called the {\it positive cone} of the representation~$U$. 
\end{itemize}

\section{Endomorphisms of standard subspaces} 
\mlabel{sec:2}

For a standard subspace 
$V \subeq \cH$, we are interested in the closed subsemigroup 
\[ S_V = \{ U \in \U(\cH) \: UV \subeq V \}\] 
of the unitary group. In the forthcoming sections, we shall study 
this semigroup by intersecting with finite dimensional subgroups 
of $\U(\cH)$. In the present section we discuss it on the general 
level to develop and present some tools that we shall use below. 
In Subsection~\ref{subsec:2.1} we show that $S_V$ is a group 
if and only if $\Delta_V$ is bounded, so that the situation is only 
interesting if $\Delta_V$ is unbounded. 
To understand the semigroup $S_V$,  it is natural to consider the full algebra 
$\cA_V := \{ A \in B(\cH) \: A V \subeq V\}$ of {\it $V$-real operators}, 
which contains $S_V$ as $\cA_V \cap \U(\cH)$. 
In Subsection~\ref{subsec:2.2} we recall an important characterization 
of the elements of $\cA_V$ in terms of the orbit maps 
$\alpha^A(t) := \Delta_V^{-it/2\pi} A \Delta_V^{it/2\pi}$ defined by the 
unitary group generated by $\Delta_V$: By results of Araki and Zsid\'o 
\cite{AZ05}, 
$A \in \cA_V$ is equivalent to the existence of an analytic continuation 
of $\alpha^A$ from $\R$ to the closed strip $\oline{\cS_\pi}$ satisfying 
$\alpha^A(\pi i) = J_V A J_V$. 
We thus obtain for every $z \in \oline{\cS_\pi}$ an injective 
representation $\alpha_z$ on $\cH$, and for 
$z = \frac{\pi i}{2}$ we even obtain a 
$*$-representation $\alpha_{\frac{\pi i}{2}} \: \cA_V \to B(\cH^{J_V}), 
A \mapsto \hat A$ by operators commuting with $J_V$. 
We conclude this section 
with Subsection~\ref{subsec:2.3}, where we take a brief look 
at one-parameter semigroups of contractions in~$\cA_V$. 

\subsection{The case where all unitary endomorphisms are invertible} 
\mlabel{subsec:2.1}

To understand the subsemigroups $S_V \subeq \U(\cH)$, one needs to 
understand when they are trivial in the sense that they are groups. 
This case is characterized in the following lemma which shows that 
standard subspaces with bounded modular operators $\Delta_V$ are too rigid 
to have non-trivial unitary endomorphisms. 

In the proof we shall need the {\it ``complementary'' standard subspace}
\[  V' := (iV)^{\bot_\R} 
= \{ \xi \in \cH \: (\forall v \in iV)\ \Re \la v,\xi \ra = 0\} 
= \{ \xi \in \cH \: (\forall v \in V)\ \Im\la v,\xi \ra = 0\}.\] 
Then 
\begin{equation}
  \label{eq:dualspace}
\Delta_{V'}  = \Delta_V^{-1}, \qquad J_V= J_{V'} \quad \mbox{ and } \quad 
J_V V= V'  
\end{equation} 
(\cite[Prop.~3.2]{Lo08}, \cite[Lemma~3.7]{NO17}).

\begin{lem} \mlabel{lem:3.2} 
For $V \in \Stand(\cH)$, the following are equivalent: 
\begin{itemize}
\item[\rm(a)] $\Delta_V$ is bounded. 
\item[\rm(b)] $V + i V = \cH$. 
\item[\rm(c)] If $H \supeq V$ is standard, then $H = V$. 
\item[\rm(d)] If $H \subeq V$ is standard, then $H = V$. 
\item[\rm(e)] The closed subsemigroup 
$S_V \subeq \U(\cH)$ is a group.
\end{itemize}
\end{lem} 

These conditions are in particular satisfied if $\cH$ is finite dimensional. 
Here (e) corresponds to the well-known fact that 
every closed subsemigroup of the compact group $\U_n(\C)$ is a group; 
cf.~also Proposition~\ref{prop:5.8}. 

\begin{prf} The  equivalence of (a) and (b) follows from $V + i V = \cD(\Delta_V^{1/2})$. 

\nin (b) $\Rarrow$ (c): If $\cH = V + i V\cong V \oplus i V$ and $H \supeq V$ is standard, 
then $\cH = H \oplus i H$ implies $V = H$. 

\nin (c) $\Leftrightarrow$ (d): Follows from $H \subeq V$ if and only if 
$V' \subeq H'$ and $\Delta_{V'} = \Delta_V^{-1}$. 

\nin (d) $\Rarrow$ (e): For $U \in S_V$ the relation $UV \subeq V$ implies $UV = V$ 
by (d) because $UV$ is also standard. 
Then $U^{-1}V = V$ as well, so that $U^{-1} \in S_V$. This shows that 
$S_V$ is a group. 

\nin (e) $\Rarrow$ (d): We show that, if $H \subeq V$ is a proper standard 
subspace, then $S_V$ is not a group. In fact, the unitary 
operator $U := J_H J_V$ satisfies 
$UV = J_H J_V V = J_H V' \subeq J_H H' = H \subeq V$. 
Therefore $U \in S_V$, and since $UV$ is a proper subset of $V$, 
the inverse $U^{-1}$ is not contained in~$S_V$. 

\nin (d) $\Rarrow$ (a): We show that, if $\Delta_V$ is unbounded, then $V$ contains a 
proper standard subspace~$V_1$. 

\nin{\bf Step 1:} First we show that $\cD(\Delta_V^{1/2}) \not\subeq \cD(\Delta_V^{-1/2})$. 
If this is not the case, then 
\[ J_V\cD(\Delta_V^{-1/2})  = \cD(\Delta_V^{1/2})\subeq \cD(\Delta_V^{-1/2})\]
 implies that 
$\cD(\Delta_V^{-1/2})$ is $J_V$-invariant. Since $J_V$ is an involution, this leads to 
$\cD(\Delta_V^{-1/2}) = J_V \cD(\Delta_V^{1/2}) = \cD(\Delta_V^{1/2})$, contradicting 
the unboundedness of $\Delta_V$. 

\nin {\bf Step 2:} By Step 1, there exists a non-zero $v_0 \in V \setminus \cD(\Delta_V^{-1/2})$ 
because $\cD(\Delta_V^{1/2}) = V + i V$. We consider the closed real hyperplane 
\[ V_1 := \{ w \in V \: \Re \la w,v_0\ra = 0\} = v_0^{\bot_\R} \cap V \subeq V.\] 
Then $V_1 \cap i V_1 \subeq V \cap i V = \{0\}$. Further, the subspace  
$V_1^{\bot_\R} = V^{\bot_\R} \oplus \R v_0 = i V' \oplus \R v_0$ 
is a real form of 
$V' + i V' + \C v_0 = \cD(\Delta_V^{-1/2}) \oplus \C v_0,$ 
so that 
$(V_1 + i V_1)^{\bot_\R} 
= V_1^{\bot_\R} \cap i V_1^{\bot_\R} = \{0\},$ 
and this implies that $V_1 + i V_1$ is dense in $\cH$. 
\end{prf}

Since every standard subspace $V$ is uniquely determined 
by the pair $(\Delta_V, J_V)$, we have: 

\begin{lem} \mlabel{lem:1.1} For $V \in \Stand(\cH)$, the stabilizer 
in the unitary group coincides with the centralizer of the pair 
$(\Delta_V, J_V)$: 
\[ \U(\cH)_V := \{ U \in \U(\cH) \: UV = V\} 
= \{ U \in \U(\cH) \: U J_V U^{-1} = J_V, U \Delta_V U^{-1} = \Delta_V \}.\] 
\end{lem}

\subsection{The algebra $\cA_V$ of $V$-real operators} 
\mlabel{subsec:2.2}

Although we are primarily interested in the subsemigroup 
$S_V \subeq \U(\cH)$, it is of some advantage to 
consider also the closed real subalgebra 
\[ \cA_V = \{ A \in B(\cH) \: A V \subeq V\} \] 
of {\it $V$-real operators}. 
The following characterization of the elements of 
$\cA$ in terms of analytic continuation of orbits maps 
(\cite[Thm.~3.18]{Lo08}, \cite[Thm.~2.12]{AZ05}) 
will be a central tool in the following. A proof can be found in 
Appendix~\ref{app:a} (Theorem~\ref{thm:lo3.18-app}). 

\begin{thm} \mlabel{thm:lo3.18} 
{\rm(Araki--Zsid\'o Theorem on $V$-real operators)} 
For $A \in B(\cH)$, the following are equivalent: 
  \begin{itemize} 
  \item[\rm(i)] $A \in \cA_V$, i.e., $A V \subeq V$. 
  \item[\rm(ii)] $A^\sharp := J_VA^*J_V \in \cA_V$. 
  \item[\rm(iii)] $\Delta_V^{1/2} A \Delta_V^{-1/2}$ is defined on 
$\cD(\Delta_V^{-1/2})$ and coincides there with $J_V A J_V$. 
  \item[\rm(iv)] The map $\alpha^A \: \R \to B(\cH), 
\alpha_t(A) := \alpha^A(t) 
:= \Delta_V^{-it/2\pi} A \Delta_V^{it/2\pi}$ extends to a strongly 
continuous function on the closed strip 
$\oline{\cS_\pi} = \{ z \in \C \: 0 \leq \Im z \leq \pi\}$ such that 
$\alpha^A$ is holomorphic on $\cS_\pi$ and $\alpha^A(\pi i) = J_VAJ_V$. 
  \end{itemize}
If these conditions are satisfied, then 
\begin{itemize}
\item[\rm(a)] $\|\alpha^A(z)\| \leq \|A\|$ for $z \in \oline{\cS_\pi}$
\item[\rm(b)] $\alpha^A(z + t) = \alpha_t(\alpha^A(z)) = 
\Delta_V^{-it/2\pi} \alpha^A(z) \Delta_V^{it/2\pi}$ 
for $z \in \oline{\cS_\pi}, t \in \R$. 
\item[\rm(c)] $\alpha^A(\oline z + \pi i) = J_V  \alpha^A(z) J_V$ 
for $z \in \oline{\cS_\pi}$. 
\item[\rm(d)] $\alpha^A(t)V \subeq V$ and 
$\alpha^A(t + \pi i)V' \subeq V'$ for all $t \in \R$. 
\end{itemize}
\end{thm}

Based on the Araki--Zsid\'o Theorem, we obtain the following remarkable 
fact, which characterizes in particular 
invertible elements in $S_V$ as those commuting either with $J_V$ or 
with~$\Delta_V$. 

\begin{cor} \mlabel{cor:degree0}
For a standard subspace $V \in \Stand(\cH)$ and $A \in B(\cH)$, 
the following are equivalent: 
\begin{itemize}
\item[\rm(i)] $AV \subeq V$ and $A$ commutes with $(\Delta_V^{it})_{t \in \R}$.
\item[\rm(ii)] $AV \subeq V$ and $A$ commutes with $J_V$. 
\item[\rm(iii)] $A$ commutes with $J_V$ and $(\Delta_V^{it})_{t \in \R}$.
\end{itemize}
It follows in particular that 
\begin{equation}
  \label{eq:unitgroup1}
 \U(\cH)_V = \{ g \in S_V \: g\Delta_V g^{-1} = \Delta_V \} 
= \{ g \in S_V \: gJ_V  g^{-1} = J_V \}. 
\end{equation}
\end{cor}

\begin{prf} (i) $\Rarrow$ (ii): If (i) is satisfied, then 
the function $\alpha^A$ is constant. Hence 
$A = \alpha^A(i\pi) = J_V \alpha^A(0) J_V = J_V A J_V$ implies (ii). 

\nin (ii) $\Rarrow$ (iii): If (ii) holds, then 
$\alpha^A(i\pi) = J_V \alpha^A(0) J_V = J_V A J_V = A = \alpha^A(0)$, 
so that Theorem~\ref{thm:lo3.18}(b) implies that 
$\alpha^A(t + \pi i) = \alpha^A(t)$ for all $t \in \R$. 
Therefore $\alpha^A$ extends to a $\pi i$-periodic bounded 
holomorphic function on all of $\C$. Now Liouville's Theorem 
implies that $\alpha^A$ is constant, so that 
$A$ commutes with $(\Delta_V^{it})_{t \in \R}$. 

\nin (iii) $\Rarrow$ (i): Condition (iii) implies that the constant 
map $\alpha^A(z) = A$ satisfies all requirements of 
Theorem~\ref{thm:lo3.18}(iv), so that $A \in \cA_V$. 

Finally, \eqref{eq:unitgroup1} follows from the equivalence 
of (i) and (ii) and Lemma~\ref{lem:1.1}. 
\end{prf}

\begin{cor} \mlabel{cor:sv1} 
The semigroup $S_V$ is invariant under the involution $\sharp$, so that 
$(S_V, \sharp)$ is an involutive semigroup. Its unitary group 
\[ \U(S_V,\sharp) := \{ s \in S_V \: s^\sharp s = s s^\sharp = \1\}  
= \{ s \in S_V \:  s^\sharp = s^{-1}\} 
= \{ s \in S_V \:  J_V s J_V = s\} 
\] 
coincides with its unit group $\U(\cH)_V = S_V \cap S_V^{-1}$. 
\end{cor}

\begin{prf} That $S_V$ is $\sharp$-invariant follows from 
Theorem~\ref{thm:lo3.18}(ii). 
Clearly, $\U(S_V,\sharp)$ consists of units of~$S_V$. 
Conversely, any $U \in S_V \cap S_V^{-1}$ satisfies $UV = V$, so that 
$U J_V U^{-1} = J_V$ and thus $U^\sharp = U^{-1}$. 
\end{prf}

The following proposition is a key tool in the following. 
It provides an analytic interpolation between the representation 
$U$ of $S_V$ on $V$ by isometries 
and an an involutive $*$-representation 
by contractions on the real subspace $\cH^{J_V}$. 

\begin{prop} \mlabel{prop:2.8} 
For every $z \in \oline{\cS_\pi}$, the map 
\[ \alpha_z \:  (\cA_V,\sharp) \to (B(\cH),\sharp), \quad A \mapsto \alpha^A(z) \] 
is an injective contractive 
morphism of real involutive unital Banach algebras with the following 
properties: 
\begin{itemize}
\item[\rm(i)] The restriction of $\alpha_z$ 
to the closed unit ball $B_V := \{ A \in \cA_V \: \|A\| \leq 1\}$ 
is continuous with respect to the strong operator topology on $B_V$  
and $B(\cH)$. 
\item[\rm(ii)] For $z = \pi i/2$, we have 
$\alpha_{\frac{\pi i}{2}}(A)^* =\alpha_{\frac{\pi i}{2}}(A^\sharp)$ 
and $\alpha_{\frac{\pi i}{2}}(A) J_V = J_V \alpha_{\frac{\pi i}{2}}(A)$, so that $\alpha_{\frac{\pi i}{2}}$ 
defines a $*$-representation of $\cA_V$ on the real Hilbert space $\cH^{J_V}$. 
\end{itemize}
\end{prop}

\begin{prf} Clearly, $\alpha_t$ is multiplicative for every $t \in \R$, 
so that  
\[ \alpha^{AB}(t) = \alpha^A(t) \alpha^B(t) \quad \mbox{ for } \quad t \in \R, 
A, B \in \cA_V.\] 
For $\xi, \eta \in \cH$, the maps 
\[ z \mapsto \la \xi, \alpha^{AB}(z) \eta \ra \quad \mbox{ and } \quad  
  z \mapsto \la \xi, \alpha^A(z) \alpha^B(z)  \eta \ra \] 
are continuous because $\alpha^A$ and $\alpha^B$ are strongly continuous 
and bounded. As both functions are holomorphic on $\cS_\pi$ 
and coincide on $\R$, 
they coincide on $\oline{\cS_\pi}$ for all $\xi,\eta\in \cH$. 
This implies that $\alpha^{AB}(z) = \alpha^A(z) \alpha^B(z)$ for all 
$z \in \oline{\cS_\pi}$.

For $A \in \cA_V$, we have 
\[ \alpha^{A^\sharp}(t) 
= \Delta_V^{-it/2\pi} J_V A^* J_V \Delta_V^{it/2\pi} 
= J_V \Delta_V^{-it/2\pi} A^* \Delta_V^{it/2\pi} J_V 
= J_V \alpha^{A}(t)^* J_V = \alpha^{A}(t)^\sharp,\] 
and therefore $\alpha^{A^\sharp}(z) = \alpha^{A}(z)^\sharp$ for $z \in \oline{\cS_\pi}.$

Now we show that $\alpha_z$ is injective. If $\alpha_z(A)= \alpha^A(z) = 0$, 
then $\alpha^A(z + t)= 0$ for all $t \in \R$, and by 
analytic continuation we get $\alpha^A=0$. 
In particular, $A =\alpha^A(0) = 0$. 

\nin (i) We have to show that, 
for $\xi \in \cH$, the map 
\[ \gamma_\xi \: B_V  \to \cH, \quad \gamma_\xi(A) := \alpha_z(A)\xi \] 
is continuous with respect to the strong operator topology on $\cA_V$.
As the linear map \break 
$\cH \to \ell^\infty(B_V,\cH), \xi \mapsto \gamma_\xi$ 
satisfies $\|\gamma_\xi\|_\infty \leq \|\xi\|$ 
(Theorem~\ref{thm:lo3.18}(a)), 
it suffices to assume that 
$\xi$ has finite spectral support with respect to the 
selfadjoint operator $\log(\Delta_V)$. Then 
\[ \alpha_z(A)\xi = \Delta_V^{-iz/2\pi} A \Delta_V^{-1/2}\eta 
\quad \mbox{ for } \quad \eta := \Delta_V^{i(z-\pi i)/2\pi} \xi.\] 
By Corollary~\ref{cor:unb-cont}, the continuity of 
$\gamma_\xi$ on $B_V$ follows from the continuity of the maps  
\[ B_V \to \cH, \quad A  \mapsto 
\Delta_V^{1/2} A \Delta_V^{-1/2}\eta  = J_V  A J_V \eta 
\quad \mbox{ and } \quad A \mapsto A \Delta_V^{-1/2} \eta \] 
(Theorem~\ref{thm:lo3.18}(iii)). 

\nin (ii) For $z = \pi i/2$, $\alpha_z(A)$ commutes with ${J_V}$  
(Theorem~\ref{thm:lo3.18}(c)), and thus 
 $\alpha_z(A)^* = \alpha_z(A)^\sharp = \alpha_z(A^\sharp)$. 
\end{prf}

On $B(\cH)$ we consider the $W^*$-dynamical system defined by 
 \[ \alpha_t(A)  = \alpha^A(t)= \Delta_V^{-it/2\pi} A \Delta_V^{it/2\pi} \quad 
\mbox{ for } \quad A \in B(\cH), t \in \R.\] 
Then, for each $A \in \cA_V$, the operators 
$\alpha^A(z)$, $z \in \cS_\pi$, belong to the space $B(\cH)^\omega$ of 
$\alpha$-analytic vectors. In particular, $\hat A \in B(\cH^{J_V})$. 
Conversely, we have: 

\begin{lem} \mlabel{lem:hatconverse} 
For an operator $B \in B(\cH)$ commuting with $J_V$, 
there exists a (unique) $A \in \cA_V$ with $\hat A = B$ if and only if 
$B$ is an $\alpha$-analytic vector whose orbit map 
$\alpha^B$ extends to a holomorphic function 
on the strip $\cS_{-\pi/2, \pi/2}$ which 
extends to a strongly continuous function on the closure. 
Then $A = \alpha^B(-\pi i/2)$. 
\end{lem}

\begin{prf} If $B = \hat A = \alpha^A(\pi i/2)$, then 
$\alpha^B(z) := \alpha^A(z + \pi i/2)$ defines a holomorphic function on 
$\cS_{-\pi/2,\pi/2}$ which is strongly continuous on the closure 
and extends the orbit map of $B$. 

Suppose, conversely, that such a function $\alpha^B$ exists 
on $\oline{\cS_{-\pi/2,\pi/2}}$. 
Then the relation $J_V BJ_V = B$ implies that 
\[ J_V \alpha^B(z) J_V = \alpha^B(\oline z) \quad \mbox{ for }\quad 
|\Im z| \leq \frac{\pi}{2}, \] 
so that $\alpha^A(z) := \alpha^B(z-\pi i/2)$ defines a holomorphic function 
on $\cS_\pi$, strongly continuous on the closure, extending the 
orbit map of $A$, and which satisfies 
\[ \alpha^A(\pi i) = \alpha^B(\pi i/2) 
= J_V \alpha^B(-\pi i/2) J_V 
= J_V \alpha^A(0) J_V = J_V A J_V.\qedhere\]
\end{prf}

In the following we shall mainly work with the 
characterization of elements $A \in S_V$ in terms of the 
analytic continuation of $\alpha^A$ to $\cS_\pi$,  but the preceding 
lemma provides a second perspective: We may also 
get information on the contraction semigroup 
$\hat{S_V} \subeq B(\cH^{J_V})$ and then obtain 
elements of $S_V$ by extending for $B \in \hat{S_V}$ 
the orbit map $\alpha^B$ to $-\frac{\pi i}{2}$. 
For a contraction $B$ on $\cH^{J_V}$, 
the regularity condition of being injective with dense range 
comes naturally into play. In this regard, we 
record the following lemma. 

\begin{lem} \mlabel{lem:regops} 
Let $\cH$ be a real or complex Hilbert space. 
Then the subset $B(\cH)_{\rm reg}\subeq B(\cH)$ of 
injective operators with dense range is a multiplicative $*$-subsemigroup 
of $B(\cH)$. It consists of those operators 
$C \: \cH \to \cH$ for which the partial isometry 
$U$ in the polar decomposition $C = U e^B$, $B = B^*$ bounded from above, 
is unitary. 
\end{lem}

\begin{prf} First we observe that the injective operators and the 
operators with dense range are multiplicative subsemigroups of $B(\cH)$. 
Hence their intersection $B(\cH)_{\rm reg}$ also is a subsemigroup. 
As $C(\cH)^\bot = \ker(C^*)$ and $C^*(\cH)^\bot = \ker(C)$, 
this semigroup is $*$-invariant. 

If $C = UP$ is the polar decomposition of $C$, then 
$P = P^* \geq 0$, and $U$ is a partial isometry from $\ker(C)^\bot$ onto 
$\oline{C(\cH)}$. Therefore the operator $C$ is injective with dense range 
if and only if $U$ is unitary. Then the positive bounded operator 
$P$ is injective with dense range, so that it can be written 
as $P = e^B$ for the operator $B := \log P$ which is bounded from above. 
\end{prf}

Although all strongly continuous one-parameter semigroups 
of $B(\cH)$ which are either symmetric or unitary 
are contained in $B(\cH)_{\rm reg}$, this is not true in general, as 
the following simple example shows: 

\begin{ex} (cf.\ \cite[Ex.~II.4.31]{EN00}) 
On the Hilbert space $\cH = L^2([0,1])$ we obtain by 
\[ (U_t f)(x) :=
\begin{cases}
  f(x+t) & \text{ for } x + t \leq 1, \\
  0  & \text{ for } x + t > 1, 
\end{cases}\] 
a strongly continuous contraction semigroup for which 
all operators $U_t$, $t > 0$, are nilpotent. For $Nt > 0$ we have 
$(U_t)^N = U_{tN} = 0$. 
\end{ex}

\begin{probs} \mlabel{prob:reg} Let $V \subeq \cH$ be a standard subspace. 

\nin (a) Show that every one-parameter semigroup 
$(U_t)_{t \geq 0}$ of $S_V$ satisfies $\hat{U_t} \in B(\cH^{J_V})_{\rm reg}$ 
for $t\geq 0$ or find an example where this is not the case. 

\nin (b) Let $B = B^* 
= e^{-H} \in B(\cH^{J_V})$ be a regular positive contraction 
for which a unitary $A \in S_V$ with $\hat A = B$ exists. 
Then the same is true for all powers 
$B^n = \hat{A^n}$, $n \in \N$, but what about the other operators 
$B^t = e^{-tH}$ for $t \geq 0$? Are they also contained in $\hat{S_V}$? 
See also Example~\ref{ex:5.7} for related problems.
\end{probs}

\subsection{One-parameter semigroups in $\cA_V$} 
\mlabel{subsec:2.3}

Classically, bounded strongly continuous one-parameter semigroups 
on Banach spaces are studied through their infinitesimal generators 
and their resolvents. We start our analysis in this subsection by 
recalling some key facts on one-parameter semigroups 
from~\cite{EN00}. This provides some tools used below 
for one-parameter subsemigroups of finite dimensional semigroups.

\begin{rem} \mlabel{rem:contr-semigroup}
(a) If $(U_t)_{t \geq 0}$ is a strongly continuous 
one-parameter semigroup of contractions on the Banach space $X$ 
and $A \: \cD(A) \to X$ its infinitesimal generator, then we have 
for every $\lambda \in \C$ with $\Re \lambda > 0$ an integral 
formula for the resolvent:  
\begin{equation}
  \label{eq:resolv}
 R(\lambda,A) := (\lambda \1 - A)^{-1} = \int_0^\infty e^{-t\lambda} U_t\, dt 
\quad \mbox{ and }\quad 
\|R(\lambda,A)\| \leq \frac{1}{\Re \lambda}
\end{equation}
(\cite[Thm.~II.1.10]{EN00})

\nin (b) If, conversely, $A \: \cD(A) \to X$ is a closed, densely defined operator on 
$X$ such that, for $\lambda > 0$,  the operators $\lambda\1 - A \: \cD(A) \to X$ have bounded 
inverses $R(\lambda,A)$ satisfying $\|R(\lambda,A)\| \leq \lambda^{-1}$, 
then $A$ is the infinitesimal generator of a uniquely 
determined semigroup of contractions (\cite[Thm.~II.3.5]{EN00}). 
That this semigroup can actually be obtained as the strong limit 
\begin{equation}
  \label{eq:euler}
U_t = \lim_{n \to \infty} \Big(\1 - \frac{t}{n}A\Big)^{-n} \quad \mbox{ for } \quad 
t > 0
\end{equation} 
follows from the discussion in \cite[\S 12.3]{HP57} 
(see also \cite{Ch68} for related results). Note that 
our assumption on $A$ implies that 
\[\Big\|\Big(\1 - \frac{t}{n}A\Big)^{-1}\Big\| 
= \frac{n}{t} \Big\|\Big(\frac{n}{t}\1 - A\Big)^{-1}\Big\| 
\leq 1,\] 
so that the right hand side of \eqref{eq:euler} is a contraction whenever 
the limit exists. 

\nin (c) If $X$ is a Hilbert space and $A$ a normal operator, then 
the assumption 
on $A$ implies that
 $\Spec(A) \subeq \C_\ell := \{ z \in \C \: \Re z \leq 0\}$. 
Then, for any $z \in \C_\ell$, 
we have $(1-tz/n)^{-n} \to e^{tz}$ for $t \geq 0$, 
as a pointwise limit of bounded functions on $\C_\ell$. Therefore 
\eqref{eq:euler} is an immediate consequence of the measurable spectral 
calculus and a normal operator $A$ generates a one-parameter semigroup 
of contractions if and only if $\Spec(A) \subeq \C_\ell$. 

\nin (d) A linear operator $A \: \cD(A) \to X$ on a Banach space 
is said to be {\it dissipative} if 
\[ \|(\lambda\1 - A)\xi\| \geq \lambda \|\xi\| 
\quad \mbox{ for } \quad \lambda > 0, \xi \in \cD(A),\] 
which is equivalent to 
\begin{equation}
  \label{eq:dissi2}
\|(\1 - hA)\xi\| \geq \|\xi\| 
\quad \mbox{ for } \quad h > 0, \xi \in \cD(A).
\end{equation}
According to the Lumer--Phillips Theorem (\cite[Thm.~II.3.15]{EN00}), 
a closed densely defined operator $A$ generates a contraction 
semigroup if and only if it is dissipative and 
$\lambda\1 - A$ has dense range for some (hence for all) 
$\lambda > 0$. If $X$ is a Hilbert space, then \eqref{eq:dissi2} implies that 
$A$ is dissipative if and only if 
\[ \Re\la \xi, A \xi \ra \leq 0 \quad \mbox{ for all } \quad \xi \in \cD(A) \] 
(\cite[Prop.~II.3.23]{EN00}).  
\end{rem}

For a standard subspace $V \in \Stand(\cH)$, we recall the 
subalgebra $\cA_V \subeq B(\cH)$ from Theorem~\ref{thm:lo3.18}. 
We are mainly interested in the semigroup $S_V = \cA_V \cap \U(\cH)$, and 
for that the representations $\alpha_z$ of $\cA_V$ will be extremely helpful. 
As $\cA_V$ is strongly closed, \eqref{eq:resolv} and \eqref{eq:euler} in 
Remark~\ref{rem:contr-semigroup} lead to: 

\begin{prop} \mlabel{prop:2.13}
{\rm(One-parameter semigroups of contractions in $\cA_V$)} 
Let $(U_t)_{t \geq 0}$ be a strongly continuous one-parameter semigroup 
of contractions on $\cH$ with  infinitesimal generator~$A$. Then 
\[ (\forall t > 0)\ \  U_t V \subeq V \qquad \Longleftrightarrow \qquad 
(\forall \lambda > 0)\ \ (\lambda \1 - A)^{-1} V \subeq V.\] 
\end{prop}

\begin{cor} \mlabel{cor:2.9} 
Let $(U_t)_{t \geq 0}$ be a strongly continuous one-parameter semigroup 
of contractions in $\cA_V$ with infinitesimal generator~$A$.
Then, for every $z \in \oline{\cS_\pi}$, 
$(\alpha_z(U_t))_{t \geq 0}$ 
is a strongly continuous one-parameter semigroup of contractions on $\cH$. 
Its infinitesimal generator $A_z$ satisfies 
\[ (\lambda \1- A_z)^{-1} = \alpha_z((\lambda \1- A)^{-1}) \quad \mbox{ for } \quad 
\lambda > 0.\] 
\end{cor}

\begin{prf} The first assertion follows immediately from 
Proposition~\ref{prop:2.8}. By Proposition~\ref{prop:2.13}, 
$(\lambda \1 - A)^{-1} \in \cA_V$ for every $\lambda > 0$, 
so that $\alpha_z((\lambda\1 - A)^{-1})$ is defined for 
$z \in \oline{\cS_\pi}$. For the second assertion we now use 
\eqref{eq:resolv} in Remark~\ref{rem:contr-semigroup} and the 
continuity of the representations $\alpha_z$. 
\end{prf}

In the following, Corollary~\ref{cor:2.9} is of particular interest for 
$z = \pi i/2$. If $A^\sharp = A$, then 
it leads to the infinitesimal generator $\hat A = A_{\pi i/2}$ 
of a symmetric contraction semigroup on $\cH^{J_V}$, showing that 
$\hat A\leq 0$. This will be important in the proof 
of the Germ Theorem (Theorem~\ref{thm:keylem3}).

The following observation will not be used below, 
but we record it here because it adds interesting information 
on certain results obtained in \cite{Ne18}, 
where we have seen that $\Stand(\cH)$ 
carries the structure of a reflection space, specified by 
$(-1)_V(W) = J_V W'$. Accordingly, a curve 
$\gamma \: \R \to \Stand(\cH)$ is called a {\it geodesics} 
if it is a morphism of reflection  spaces, where 
$\R$ carries the canonical reflection structure 
given by the point reflections $(-1)_x(y) = 2y - x$. 
By \cite[Prop.~2.9]{Ne18}, geodesics 
$\gamma \: \R \to \Stand(\cH)$ 
with $\gamma(0) = V$ for which the curve $(J_{\gamma(t)})_{t \in \R}$ 
is strongly continuous, are the curves of the form 
$\gamma(t) = U_t V,$ 
where $(U_t)_{t \in \R}$ is a unitary one-parameter group satisfying 
$J_V U_t J_V = U_{-t}$ for $t \in \R$.

\begin{prop} Assigning to the generator $A = A^\sharp = - A^*$ of a strongly 
continuous $\sharp$-symmetric unitary one-parameter semigroup in $\cA_V$ 
the curve $(e^{tA}V)_{t \in \R}$ in $\Stand(\cH)$, we obtain a bijection 
onto the set of decreasing geodesics $\gamma \: \R \to \Stand(\cH)$ 
with $\gamma(0) = V$. 
\end{prop}

\begin{prf} The relation $A^\sharp = A$ is equivalent to $A^* = {J_V}A{J_V}$. 
If, in addition, $A^* = - A$, then ${J_V}A{J_V} = - A$. Then the curve 
$\gamma(t) := e^{tA}V$ defines a geodesic in $\Stand(\cH)$  
which is decreasing because $t < s$ implies that 
$\gamma(s) = e^{tA} e^{(s-t)A} V \subeq e^{tA} V = \gamma(t)$. 
That all decreasing geodesics are of this form 
follows from \cite[Prop.~2.9]{Ne18}. 
\end{prf}

\section{Wick rotations of tubes and Olshanski semigroups} 
\mlabel{sec:3}

To apply tools from finite dimensional Lie theory,
we consider subsemigroups of $B(\cH)$ that arise by analytic 
continuation of a unitary representation $U \: G \to U(\cH)$ 
of a Lie group $G$ to a semigroup $S^U = G \exp(iC_U)$, 
where $C_U := \{ x \in \g \: -i\partial U(x) \geq 0\}$ 
is the positive cone of~$U$.
Assuming that $U$ has discrete kernel,  the semigroup $S^U$ always 
exists and $U(g \exp(ix)) = U(g) e^{i \partial U(x)}$ provides an extension 
of $U$ to $S^U$. To implement $J_V$ as well, we also  consider 
an involution $\tau_G \in \Aut(G)$ (inducing an involution $\tau$ on $\g$),  
for which $U$ extends by 
to an antiunitary representation of the graded Lie group $G \rtimes \{\id_G,\tau_G\}$. 
Then $J := U(\tau_G)$ is a conjugation satisfying 
$U(\tau_G(g)) = J U(g) J$ for $g \in G$. 
For $h \in \g^\tau$ we now consider the 
standard subspace $V$ determined by 
\[ J_V = J \quad \mbox{ and } \quad 
\Delta_V = e^{2\pi i \partial U(h)},\quad \mbox{  so that } \quad 
\Delta_V^{-it/2\pi} = U(\exp th).\]  
By the Araki--Zsid\'o Theorem, we are now led to 
the problem to determine the subsemigroup 
$S^U_{\rm inv}$ of those elements 
$s \in G$ for which the orbit map $\beta^s(t) = \exp(th) s \exp(-th)$ 
extends holomorphically to the closure of $\cS_\pi$, in such a way 
that $\beta_{\pi i}(s) = \tau_G(s)$. 
In Theorem~\ref{thm:4.17}, we 
show that 
\[ S^U_{\rm inv} \subeq S_V = \{ g \in G \: U(g)V \subeq V\}.\] 
To prepare this theorem, we start in Subsection~\ref{subsec:3.1} 
with a discussion of the ``abelian case'', where 
$G$ is simply a real Banach space $E$, endowed with an 
involution $\tau$ and an endomorphism $h$, and 
$W \subeq E$ is a pointed closed convex cone invariant 
under $-\tau$ and the one-parameter group $e^{\R h}$. 
In this simple situation the semigroup $(E + i W)_{\rm inv}$ 
is a closed convex cone in $E$ that 
can be determined very explicitly by elementary means. It provides an 
important blueprint for the more general non-abelian situation. 
In Subsection~\ref{subsec:3.2} we then recall 
the basic facts on Olshanski semigroups $\Gamma_G(W)$ 
for invariant cones $W \subeq \g$. They are the non-abelian generalizations 
of the tube $E + i W$. Finally, we verify the inclusion 
$S^U_{\rm inv} \subeq S_V$ in Subsection~\ref{subsec:3.3}.

\subsection{Wick rotations of tubes} 
\mlabel{subsec:3.1}

In this section we develop another tool that we shall use below 
in the context of Lie algebras. This subsection 
represents some key geometric features that can already 
be formulated in the abelian 
context. 

Let $E$ be a real Banach space 
endowed with the following data: 
\begin{itemize}
\item A continuous involution $\tau\in \GL(E)$; we write $E = E^+ \oplus E^-$, 
$E^\pm := \ker(\tau\mp \1)$ for the 
$\tau$-eigenspace decomposition.  
\item An endomorphism $h \in B(E)$, commuting with $\tau$. 
\item A closed convex cone $W \subeq E$ which is pointed, i.e., $W \cap -W = \{0\}$, and 
invariant under $-\tau$ and the one-parameter group~$e^{\R h}$.
\end{itemize}

We consider the closed convex cone  
\[ T_W := E + i W \subeq E_\C,\] 
which is obviously invariant under $e^{\R h}$ and $-\tau$, where we use the same 
notation for the complex linear extensions to $E_\C$. 
We do not assume that the cone $W$ has interior points, so that 
$W - W$ may be a proper subspace of $E$.
If $\sigma^c \: E_\C\to E_\C$ is the antilinear involution with 
fixed point set $E^c := E^+ + i E^-$, 
then $\sigma^c$ acts on $iE$ as $-\tau$, so that 
$\sigma^c(T_W) = T_W$, and 
\[T_W^{\sigma^c} =  T_W \cap E^c = E^+ + i (W \cap E^-) \] 
is the closed convex cone of $\sigma^c$-fixed points in $T_W$.
We are interested in the closed convex cone 
\[ T_{W,{\rm inv}} := \{ x \in E \: 
e^{yi h}x \in T_W  \ \mbox{ for } \  y\in [0,\pi];\  
e^{\pi i h} x = \tau(x) \}.\] 

\begin{lem}
  \mlabel{lem:newcone} 
{\rm(a)} For $x \in E$, the condition $e^{\pi i h} x = \tau(x)$ is equivalent to 
$e^{\frac{\pi i}{2} h} x \in E^c$. \\
{\rm(b)} $T_{W,{\rm inv}} = \{ x \in E \: 
e^{z h}x \in T_W  \ \mbox{ for } \  z\in \oline{\cS_\pi};\  e^{\pi i h} x = \tau(x) \}.$
\end{lem}

\begin{prf} (a) As above, let 
$\sigma^c \: E_\C \to E_\C$ denote the antilinear extension 
of $\tau$, 
so that $\Fix(\sigma^c) = E^c$. 
For $x \in E$ and $x^c := e^{\frac{\pi i}{2} h} x$, 
the condition $x^c \in E^c$ is equivalent to $\sigma^c(x^c) = x^c$, 
which is equivalent to 
\[  e^{\frac{\pi i}{2} \ad h} x 
= \sigma^c(e^{\frac{\pi i}{2} \ad h} x)
= e^{-\frac{\pi i}{2} \ad h} \sigma^c(x)
= e^{-\frac{\pi i}{2} \ad h} \tau(x).\] 
This in turn is equivalent to $e^{\pi i \ad h}x  = \tau(x)$. 

(b) follows from the fact that $T_W$ is invariant under $e^{\R h}$. 
\end{prf}

For an $h$-invariant 
real subspace $F \subeq E_\C$, we write $F_\lambda = F_\lambda(h) 
:= F \cap \ker(h - \lambda\1)$ for the $h$-eigenspaces in~$F$. 

\begin{lem} \mlabel{lem:4.13a} 
For $x  \in E$, the condition $e^{\pi i h} x = \tau(x)$ is equivalent to the 
existence of finitely many elements $x_n \in E_n(h)$ with 
$x = \sum_{n \in \Z} x_n$ and $\tau(x_n) = (-1)^n x_n$. 
\end{lem}

\begin{prf} We write $x = x_+ + x_-$ with $x_\pm \in E^\pm$. 
Then $e^{\pi i h} x = \tau(x)$ is equivalent to 
\begin{equation}
  \label{eq:even-odda}
 e^{\pi i h} x_+= x_+\quad \mbox{ and }  \quad 
e^{\pi i h} x_- = - x_-.
\end{equation}
Combining both, we see that $e^{2\pi i h} x =  x.$
The space $E_\C^{\rm fix}$ 
of fixed points of the automorphism $e^{2\pi i h} \in \GL(E_\C)$ 
carries a norm continuous 
action of the circle $\T \cong \R/\Z$, defined by $\beta_t(y) := e^{2\pi i t h}y$. 
As $h$ is bounded, $E_\C^{\rm fix}$ is a direct sum of finitely many 
$h$-eigenspaces $E_{\C,n}(h)$, $n \in \Z$. 
Accordingly, we write 
\[x = \sum_{n \in \Z} x_n \quad \mbox{ with } \quad hx_n = n x_n.\]
As $\|h\| < \infty$, only finitely many summands are non-zero. 
The antilinear involution $\sigma$ of $E_\C$, whose fixed point set is $E$,  
commutes with $h$. Therefore the $h$-eigenspaces 
are $\sigma$-invariant, and thus $\sigma(x) = x$ implies 
$\sigma(x_n) = x_n$ for every $n \in \Z$, i.e., $x_n \in E$. 
Now $e^{\pi i h}x_n = (-1)^n x_n$ and 
\eqref{eq:even-odda} imply that  $x_+$ is the sum of the $x_n$ with 
$n$ even, and $x_-$ is the sum of the $x_n$ with $n$~odd. 
As $\tau(x_\pm) = \pm x_\pm$, this 
in turn shows that $\tau(x_n) = (-1)^n x_n$ for $n \in \Z$. 

If, conversely, $x = \sum_{n \in \Z} x_n$ with $x_n \in E$ satisfying 
$hx_n = nx_n$ and $\tau(x_n) = (-1)^n x_n$, then 
the relation $e^{\pi i h} x = \tau(x)$ is obvious. 
\end{prf}

The following proposition is a key geometric ingredient of the proof of 
our Structure Theorem (Theorem~\ref{thm:main2}). 

\begin{prop} \mlabel{prop:keylema} 
The cone $T_{W,{\rm inv}}$ has the simple form 
\begin{equation}
  \label{eq:w-decomp}
 T_{W,{\rm inv}} = (W \cap E^-_1(h)) \oplus E^+_0(h) \oplus (-W \cap E^-_{-1}(h)).
\end{equation}
In particular, it is determined by the two cones $W \cap E^-_{\pm 1}(h)$. 
\end{prop}

\begin{prf} First we 
note that the cone $T_{W,{\rm inv}}$ is closed and 
invariant under $e^{\R h}$ because 
$T_W$ is invariant under $e^{\R h}$, the operators $e^{yi h}$ commute with 
$e^{\R h}$, and so does $\tau$. 

Let $x \in T_{W,{\rm inv}}$. Then Lemma~\ref{lem:4.13a} implies that 
$x = \sum_{n \in \Z} x_n$ is a finite sum 
with $x_n \in E_n(h)$ and $\tau(x_n) = (-1)^n x_n$.
We claim that $x_n = 0$ for $|n| > 1$. 
Suppose first that there exists an $n > 1$ with $x_n \not=0$ 
and assume that $n$ is maximal with this property. 
Then the invariance of $T_{W, {\rm inv}}$ 
under $e^{\R h}$ and its closedness imply that 
\begin{equation}
  \label{eq:asymp}
x_n 
= \lim_{t \to \infty} e^{-nt}\sum_{m \in \Z} e^{m t} x_m 
= \lim_{t \to \infty} e^{-nt} e^{t h} x \in T_{W, {\rm inv}}.
\end{equation}
As $n > 1$, we now obtain
\[ e^{[0,\pi]i h} x_n = e^{[0,n\pi]i} x_n \ni \pm i x_n.\] 
This leads to $\pm x_n \in W$, and since $W$ is pointed, we arrive at the 
contradiction $x_n = 0$. 
An analogous argument shows that $x_n = 0$ for $n < -1$. 
This shows that $x = x_1 + x_0 + x_{-1}$ 
with $x_0 \in E^+_0(h)$, and $x_{\pm 1} \in T_{W, {\rm inv}}$, 
obtained from \eqref{eq:asymp}, implies that 
$\pm i x_{\pm 1} = e^{\frac{\pi i}{2} h} x_{\pm 1} \in T_W = E + i W$, hence 
$x_{\pm 1} \in \pm W \cap E^-_{\pm 1}(h)$. 

Conversely, every element of the form 
$x = x_{-1} + x_0 + x_1$ with $x_0 \in E^+_0(h)$ and $x_{\pm 1} \in \pm W \cap 
E^-_{\pm 1}(h)$ 
is contained in $T_{W, {\rm inv}}$ because 
$e^{\pi i h} x = -x_{-1} + x_0 - x_1 = \tau(x)$, and 
\[ e^{yi h} x = \underbrace{x_0}_{\in E^+} + \underbrace{\cos(y)(x_{-1} + x_1)}_{\in  E^-} 
+ \underbrace{i \sin(y)(x_1 - x_{-1})}_{i E^-} \in E + iW,\] 
because $\sin(y)(x_1 - x_{-1}) \in W$ for $y \in [0,2\pi]$. 
\end{prf}

From Proposition~\ref{prop:keylema} we immediately obtain: 
\begin{cor} \mlabel{cor:keylema-lie} 
Let $\g$ be a finite dimensional real Lie algebra, endowed with 
an involution $\tau\in \Aut(\g)$ with eigenspace decomposition 
$\g = \fh \oplus \fq$, 
$\fh = \ker(\tau- \1)$ and $\fq := \ker(\tau + \1)$, 
an element $h \in \fh$, and a pointed closed convex cone $W \subeq \g$, 
invariant under $-\tau$ and $e^{\R \ad h}$. 
For $T_W := \g + i W$, the cone $T_{W,{\rm inv}}$ then has the simple form 
\begin{equation}
  \label{eq:w-decomp-lie}
 T_{W,{\rm inv}} = (W \cap \fq_1(h)) \oplus \fh_0(h) \oplus (-W \cap \fq_{-1}(h)).
\end{equation}
In particular, it is determined by the two pointed cones $W \cap \fq_{\pm 1}(h)$. 
\end{cor}

\subsection{Extensions of unitary representations to 
semigroups} 
\mlabel{subsec:3.2}

Below we shall need  non-abelian analogs of the tubes 
$T_W = E + i W$, where $E$ is replaced by a finite dimensional 
simply connected Lie group 
$G$ and $W \subeq \fg$ is an $\Ad(G)$-invariant closed convex cone. 

\begin{defn} \mlabel{def:olshanski} (Olshanski semigroups) 
Let $G$ be a $1$-connected Lie group with Lie algebra 
$\g$ and $W \subeq \g$ be a pointed $\Ad(G)$-invariant closed convex cone. 
\begin{footnote}{Then $W$ is {\it weakly elliptic} 
in the sense that $\Spec(\ad x) \subeq i \R$ holds for every $x \in W$. 
In fact, by \cite[Prop.~VII.3.4(b)]{Ne00} $W$ is weakly elliptic in 
the ideal $W-W$, 
and since $[x,\g] \subeq W-W$ holds for any $x \in W$, 
we have $\Spec(\ad x) \subeq \{0\} \cup \Spec(\ad x\res_{W-W})\subeq i \R$.}
\end{footnote}
The corresponding {\it Olshanski semigroup} $\Gamma_G(W)$ is defined as follows.
Let $G_\C$ be the $1$-connected Lie group with Lie algebra $\g_\C$ 
and let $\eta_G \: G \to G_\C$ be the canonical morphism of Lie groups 
for which $\L(\eta) \: \g \into \g_\C$ is the inclusion. 
\begin{footnote}{
In general the map $\eta_G$ is not injective, as the example $G = \tilde\SL_2(\R)$ with
$G_\C = \SL_2(\C)$ shows.}   
\end{footnote}
As $G_\C$ is simply connected, the complex conjugation 
on $\g_\C$ integrates to an antiholomorphic 
involution $\sigma \: G_\C \to G_\C$ with 
$\sigma \circ \eta_G = \eta_G$, and this implies that
$\eta_G(G)$ coincides with the subgroup $(G_\C)^\sigma$ 
of $\sigma$-fixed points in $G_\C$.
\begin{footnote}{Since $G_\C$ is simply 
connected, this subgroup is connected by Corollary~\ref{cor:fixedpoint}.
}\end{footnote}

As $W$ is weakly elliptic, 
Lawson's Theorem (\cite[Thm.~XIII.5.6]{Ne00}) asserts that
\[\Gamma'_G(W) := \Gamma_{(G_\C)^\sigma}(W) := (G_\C)^\sigma \exp(iW) \subeq G_\C \]
is a closed subsemigroup of $G_\C$ stable under the antiholomorphic
involution $s^* := \sigma(s)^{-1}$, and that the polar map
\[ (G_\C)^\sigma \times W \to \Gamma'_G(W), \quad (g,x) \mapsto g \exp(ix) \]
is a homeomorphism. 
Next we observe that $\ker \eta_G$ is a discrete subgroup of $G$ and define 
$\Gamma_G(W)$ as the simply connected covering of $\Gamma'_G(W)$ 
(\cite[Def.~XI.1.11]{Ne00}). 
Basic covering theory implies that $\Gamma_G(W)$ inherits an involution 
given by 
\[ (g \exp(ix))^* = \exp(ix)g^{-1} = g^{-1} \exp(\Ad(g)ix)\] 
 and a homeomorphic  polar map $G \times W \to \Gamma_G(W), 
(g,x) \mapsto g \exp(ix)$. 
We write $\exp \: \g + i W \to \Gamma_G(W)$ for the canonical
lift of the exponential function 
\[ \exp \: \L(\Gamma_G'(W)) = \g + i W \to \Gamma'_G(W) \subeq G_\C.\] 
For every $x \in \g + i W$, the curve $\gamma_x(t) := \exp(tx)$ 
is a continuous one-parameter semigroup of $\Gamma_G(W)$.  

If $W$ has interior points, then 
the polar map restricts to a diffeomorphism from $(G_\C)^\sigma \times W^0$ 
onto the interior $\Gamma'_G(W^0)$ of $\Gamma'_G(W)$. 
Further, $\Gamma_G(W^0)= G \exp(i W^0)$ is a complex manifold with a holomorphic 
multiplication and the exponential function 
$\g + i W^0 \to \Gamma_G(W^0)$ is holomorphic, 
whereas the involution $*$ is antiholomorphic 
(\cite[Thm.~XI.1.12]{Ne00}). 
\end{defn}

We now turn to the analytic continuation of unitary representations 
of $G$ to Olshanski semigroups $\Gamma_G(W)$. 

\begin{prop} {\rm(Holomorphic extension of unitary representations)} 
\mlabel{prop:ucprop-new}
Let $(U,\cH)$ be a unitary representation of $G$ with discrete kernel 
and consider the ideal $\fn_U := C_U -C_U$ and the corresponding 
normal integral subgroup $N_U \trile G$.\begin{footnote}
{Normal integral subgroups of $1$-connected 
Lie groups are always closed and $1$-connected by \cite[Thm.~11.1.21]{HN12}.}  
\end{footnote}
Then the following assertions hold: 
  \begin{itemize}
  \item[\rm(i)] $U$ extends by 
$U(g\exp(ix)) = U(g) e^{i \partial U(x)}$ 
to a strongly continuous contraction representation 
of the Olshanski semigroup $S^U := G \exp(i C_{U})$ 
which is holomorphic on the complex manifold $N_U \exp(i C_U^0)$. 
\item[\rm(ii)] If $J \: \cH \to \cH$ is a conjugation and 
$\tau_G \in \Aut(G)$ an involution with derivative $\tau \in \Aut(\g)$, 
satisfying $J U(g) J = U(\tau_G(g))$ for $g \in G$, 
then the involutive automorphism of $S^U$ given by 
$\tau_S(g \exp(ix)) = \tau_G(g) \exp(-i \tau(x))$ 
satisfies  $J U(s) J = U(\tau_S(s))$ for $s \in S^U$. 
\end{itemize} 
\end{prop}

\begin{prf} (i) The assumption that $\ker(U)$ is discrete implies that 
$C_U$ is pointed. That $U(g\exp(ix)) = U(g) e^{i \partial U(x)}$ 
defines a representation which is holomorphic and non-degenerate on 
$\Gamma_{N_U}(C_U^0) = N_U \exp(i C_U^0)$ follows from  \cite[Thm~XI.2.5]{Ne00}.  
Now \cite[Cor.~IV.1.18, Prop.~IV.1.28]{Ne00} 
imply that $U$ is strongly continuous on $\Gamma_{N_U}(C_U)$ 
because $U(\Gamma_{N_U}(C_U))$ is bounded. The continuity on $S^U$ now follows 
from $S^U = G \Gamma_{N_U}(C_U) = G \exp(iC_U)$, the fact that the polar map 
is a homeomorphism, and the strong continuity of the multiplication on 
the operator ball.  

\nin (ii) The relation $J U(g) J = U(\tau_G(g))$ implies 
$J \partial U(x) J = \partial U(\tau(x))$ for $x \in \g$, and therefore 
$J i\partial U(x) J = -i\partial U(\tau(x))$ implies that 
$-\tau(C_U) = C_U$. Therefore the involution $\tau_S(g\exp(ix)) 
= \tau_G(g) \exp(-i\tau(x))$ on $S^U$ is defined. As it is the unique 
continuous lift of an automorphism of $\Gamma'_G(C_U) \subeq G_\C$,  
preserving the base point $e \in S^U$, it defines an automorphism of~$S^U 
= \Gamma_G(C_U)$. For $s = g \exp(ix)$ we have 
\[ J U(s) J = J U(g) J J e^{i \partial U(x)} J 
= U(\tau_G(g)) e^{-i \partial U(\tau x)} = U(\tau_G(g)\exp -\tau(x)) 
= U(\tau_S(s)).\qedhere\]
\end{prf}

\begin{rem} \mlabel{rem:4.15} 
(a) Let $(U,\cH)$ be an antiunitary representation 
of the graded Lie group $(G,\eps_G)$ and 
$\sigma \in G_-$ be an involution. We write $\tau_G(g) = \sigma g \sigma$ for the 
corresponding involutive automorphism of $G$ and $\tau = \Ad(\sigma) \in \Aut(\g)$ 
for the corresponding involution of the Lie algebra. 
Then the  positive cone $C_U$ of $U$ is a closed convex cone satisfying 
\[ \Ad(g)C_U = \eps_G(g) C_U \quad \mbox{ for } \quad g \in G.\] 
In particular, it is invariant under $-\tau$ 
(cf.\ Proposition~\ref{prop:ucprop-new}(ii)).

(b)  The fixed point set of the involution $\tau_S$ on $S^U$ 
is the subsemigroup 
\[ S^U_G := G^{\tau_G}  \exp(i(\fq \cap C_{U})) \subeq S^U \] 
because $s = g \exp(ix)$ is $\tau_S$-invariant if and only if 
$\tau_G(g) = g$ and $\tau(x) = -x$, i.e., $x \in \fq \cap C_U$.
\end{rem}

\subsection{Wick rotations of Olshanski semigroups} 
\mlabel{subsec:3.3}

In this subsection we describe how holomorphic extensions of 
unitary representations of complex Olshanski semigroups can 
be used to obtain non-trivial endomorphism semigroups $S_V \subeq G$ for 
certain standard subspaces. 

Let $G$ be a $1$-connected Lie group and 
$W \subeq \g$ be a pointed invariant closed convex cone, 
so that we have the Olshanski semigroup $\Gamma_G(W) = G \exp(iW)$ 
which is the simply connected covering of the semigroup 
$\Gamma'_G(W)\subeq G_\C$. We write 
$q_S \: \Gamma_G(W) \to \Gamma'_G(W) \subeq G_\C$ for the 
universal covering map (Definition~\ref{def:olshanski}). 
We further assume that $\tau_G \in \Aut(G)$ is an involution and that 
the corresponding automorphism $\tau \in \Aut(\g)$ satisfies $\tau(W) = -W$. 
For an element $h \in \fh = \g^\tau$, we consider the $\R$-action on $\Gamma_G(W)$, given by 
\[  \beta_t(s) := \beta^s(t) := \exp(th) s \exp(-th)\quad \mbox{ for }  
\quad s \in \Gamma_G(W), t \in \R\]  
and note that the corresponding $\R$-action on $G_\C$ extends 
to a holomorphic $\C$-action by 
\begin{equation}
  \label{eq:betaz}
 \beta_z(g) := \beta^g(z) := \exp(zh) g \exp(-zh)\quad \mbox{ for }  
z \in \C, g \in G_\C.
\end{equation}

\begin{defn}  \mlabel{def:holodisc}
If $M$ is a complex manifold, then we call a continuous map 
$f \: M \to \Gamma_G(W)$ {\it holomorphic} if the composition 
$q_S \circ f \: M \to G_\C$ is holomorphic. 

Holomorphic extensions of the orbits maps $\beta^s \: \R \to \Gamma_G(W)$ have 
to be understood in this sense. 
For $z \in \C$, we say that $\beta^s(z)$ {\it exists} 
if there exists a closed strip $\oline{\cS_{a,b}} \subeq \C$ containing $\R$ 
and~$z$, and an extension of $\beta^s$ to a continuous map 
$\oline{\cS_{a,b}} \to \Gamma_G(W)$ which is holomorphic on $\cS_{a,b}$. 
Then we write $\beta^s(z) = \beta_z(s)$ 
for the value of this analytic continuation in~$z$, 
which does not depend on the choice of $a$ and $b$ as long as 
$a \leq \Im z \leq b$. 
\end{defn}

\begin{lem} \mlabel{lem:wick-olsh}
For $z \in \C$, let $\Gamma_G(W)_{z} \subeq \Gamma_G(W)$ be the set of all elements $s\in \Gamma_G(W)$ for 
which $\beta^s(z)$ exists. Then the following assertions hold: 
  \begin{itemize}
\item[\rm(i)] $\Gamma_G(W)_{z} = q_S^{-1}(\Gamma'_G(W)_{z})$ and $q_S \circ \beta_z = \beta_z \circ q_S$ 
on $\Gamma_G(W)_{z}$. 
  \item[\rm(ii)] $\Gamma_G(W)_z$ is a closed subsemigroup of $\Gamma_G(W)$ and 
$\beta_z \: \Gamma_G(W)_z \to  \Gamma_G(W)$ is a continuous homomorphism. 
  \item[\rm(iii)] The closed subsemigroup $\Gamma_G(W)^{\tau_S} := G^{\tau_G} \exp(i(\fq \cap W))$ 
is the set of fixed points of the involutive automorphism 
$\tau_S$ of $\Gamma_G(W)$, defined by $\tau_S(g \exp(ix)) = \tau_G(g) \exp(-i\tau(x))$. 
\item[\rm(iv)] $\Gamma_G(W)_{{\rm inv}} := 
\{ g \in G \cap \Gamma_G(W)_{\pi i} \: \beta_{\pi i}(g) = \tau_G(g)\}$ 
is a closed subsemigroup of $\Gamma_G(W)$ with Lie wedge 
\[ \L(\Gamma_G(W)_{{\rm inv}}) 
= \L(\Gamma_G(W))_{\rm inv} = (\g + i W)_{\rm inv}.\] 
  \end{itemize}
\end{lem}

We recall from Corollary~\ref{cor:keylema-lie} 
the explicit description of $(\g + i W)_{\rm inv}$ as 
\[ (\g + i W)_{\rm inv} 
= (W \cap \fq_1(h)) \oplus \fh_0(h) \oplus (-W \cap \fq_{-1}(h)).\] 

\begin{prf} (i) Since $q_S \circ \beta_t = \beta_t \circ q_S$ holds for $t \in \R$, 
the uniqueness of analytic continuation implies that 
$q_S(\Gamma_G(W)_{z}) \subeq \Gamma'_G(W)_z$ with 
\[ q_S \circ \beta_z = \beta_z \circ q_S \: \Gamma_G(W)_{z} \to G_\C.\] 
If $s \in \Gamma_G(W)$ is such that $q_S(s) \in 
\Gamma'_G(W)_{z}$, we fix an analytic 
continuation 
$\beta^{q_S(s)} \: \oline{\cS_{a,b}} \to \Gamma'_G(W) \subeq G_\C$ 
of the orbit map $\beta^{q_S(s)} \: \R \to \Gamma'_G(W)$. 
As the closed strip $\oline{\cS_{a,b}}$ is simply connected, there exists a 
unique continuous lift 
\[ \tilde\beta^s \: \oline{\cS_{a,b}} \to \Gamma_G(W) \quad \mbox{ with } \quad 
q_S \circ \tilde\beta^s = \beta^{q_S(s)} \circ q_S \quad \mbox{ and } \quad 
\tilde\beta^s(0) = s.\] 
Then the uniqueness of continuous lifts to coverings implies that 
$\tilde\beta^s(t) = \beta^s(t)$ for $t \in \R$ and, 
by construction, $\tilde\beta^s$ is holomorphic on $\cS_{a,b}$. 
This implies that $s \in \Gamma_G(W)_{z}$ with $\beta^s(z) = \tilde\beta^s(z)$. 

\nin (ii) In $G_\C$ we have for $\Im z > 0$ that 
\[ \Gamma'_G(W)_{z} 
= \{ s \in \Gamma'_G(W) \: 0 \leq \Im w \leq \Im z \Rarrow \beta_w(s) \in \Gamma'_G(W) \} 
= \bigcap_{0 \leq \Im w \leq \Im z} \beta_w^{-1}(\Gamma'_G(W)),\] 
where $\beta_w \in \Aut(G_\C)$ is the unique automorphism 
from \eqref{eq:betaz} with 
$\L(\beta_w) = e^{w \ad h}$. Since $\Gamma'_G(W)$ is a closed subset of $G_\C$, 
the subset $\Gamma'_G(W)_{z}$ of $\Gamma'_G(W)$ is closed. Now (i) implies that 
$\Gamma_G(W)_z = q_S^{-1}(\Gamma'_G(W)_z)$ is also closed. 

Next we show that $\Gamma_G(W)_{z}$ is a subsemigroup on which $\beta_z$ is multiplicative. 
Let $s_1, s_2 \in \Gamma_G(W)_{z}$ and consider the minimal strip 
$\cS_{a,b} \subeq \C$ with $\R \cup \{z\} \subeq \oline{\cS_{a,b}}$. 
Then the map 
\[ \beta^{s_1} \cdot \beta^{s_2} \: \oline{\cS_{a,b}} \to \Gamma_G(W), 
\quad z \mapsto \beta^{s_1}(z) \beta^{s_2}(z), \] 
is continuous and holomorphic on $\cS_{a,b}$. For $t \in \R$, we have
$(\beta^{s_1} \beta^{s_2})(t) 
= \beta_t(s_1) \beta_t(s_2) = \beta_t(s_1 s_2)$ 
because $\beta_t$ is an automorphism of $\Gamma_G(W)$. 
Uniqueness of analytic continuation therefore implies that 
$s_1 s_2 \in \Gamma_G(W)_{z}$ with $\beta_z(s_1 s_2) = \beta^{s_1}(z) \beta^{s_2}(z) 
= \beta_z(s_1) \beta_z(s_2)$. 

\nin (iii) On $G_\C$ we have a unique antiholomorphic involution $\sigma^c$ 
inducing on $\g_\C$ the antilinear extension of $\tau$, so that 
its group of fixed points has the Lie algebra $\g^c$.
It acts on $s = g \exp(ix)$ by 
$\sigma^c(s) = \tau_G(g) \exp(-i \tau(x))$. 
By uniqueness of lifts to coverings, 
this implies that $\tau_S$ defines an involutive automorphism of 
$\Gamma_G(W)$, and the assertion 
follows immediately from the formula for $\tau_S$. 
  
\nin (iv) follows from the trivial observation 
that, for a family $(S_j)_{j \in J}$ of closed subsemigroups 
of $\Gamma_G(W)$, we have $\L\big(\bigcap_j S_j\big) = \bigcap_j \L(S_j)$. 
\end{prf}

The following lemma provides an interesting tool that permits 
us to work effectively with holomorphic maps with values 
in $\Gamma_G(W)$, which neither is a manifold nor ``complex''. 

\begin{lem} \mlabel{lem:autoholo} 
Let $W \subeq \g$ be a closed pointed 
convex invariant cone and 
\[ U \: \Gamma_G(W) \to B(\cH), \quad U(g\exp(ix)) = U(g) e^{i\partial U(x)} \] 
be a $*$-representation obtained from a unitary representation $U$ of $G$. 
Then, for every holomorphic map $f \: M \to \Gamma_G(W)$, 
$M$ a finite dimensional complex manifold, 
the composition \break $U \circ f \: M \to B(\cH)$ is holomorphic. 
\end{lem}

\begin{prf} As the assertion is local with respect to $M$, we 
may w.l.o.g.\ assume that $M$ is connected and that $f(M)$ has compact closure 
in $\Gamma_G(W)$. 
Let $f \: M \to \Gamma_G(W)$ be a holomorphic map. 
Then $f' := q_S\circ f \:  M \to \Gamma'_G(W) = G_\C$ is holomorphic 
by definition. 
We consider the ideal $\fn := W - W \trile \g$ and the corresponding 
normal integral subgroup $N \trile G$. As $N$ is closed and 
$1$-connected by \cite[Thm.~11.1.21]{HN12}, we obtain a 
quotient group $Q := G/N$. We likewise have a closed normal subgroup 
$N_\C \trile G_\C$ and  $Q_\C := G_\C/N_\C$. 
Let $r : G_\C \to Q_\C$ denote the quotient map. 
Then 
\[ r(\Gamma_G'(W)) = r(G \exp(iW)) = r(G) \subeq (Q_\C)^\sigma \subeq Q_\C \] 
is contained in the totally real submanifold $(Q_\C)^\sigma$ 
of fixed points of the 
antiholomorphic involution $\sigma$ of $Q^\C$ corresponding to the 
complex conjugation of $\fq_\C$ with respect to~$\fq$. 
We conclude that the holomorphic map $r \circ f' \: M \to r(G) \subeq Q_\C$ 
is constant, hence equal to $r(g')$ for some $g' \in q_S(G) \cap f'(M)$. 
This implies the existence of a holomorphic map 
$h \: M \to N_\C$ with 
$f'(m) = g' h(m)$ for $m \in M$. Lifting to the covering 
space $\Gamma_G(W)$, we conclude that 
there exists a $g \in G$ with $f(m) = g h(m)$ for all $m \in M$, 
where $h \: M \to \Gamma_N(W)$ is a holomorphic map. 

Pick $x \in W^0$ (the interior of $W$ with respect to $\fn$) 
and put $s_n := \exp(n^{-1}x) \in \Gamma_G(W^0)$. We thus obtain a 
sequence $U(s_n)^* = U(s_n)$ of hermitian operators converging strongly to~$\1$.
Further, $s_n h(M)$ is contained in the complex manifold $\Gamma_N(W^0)$ 
and the map $h_n \: M \to \Gamma_N(W^0), m \mapsto s_n h(m)$ 
is holomorphic. Therefore the maps $H_n := U \circ h_n \: M \to B(\cH)$ 
are holomorphic and we want to show that $H := U \circ h$ is also holomorphic. 
For $\xi,\eta \in \cH$, we have 
\[ \lim_{n \to \infty} \la \xi, H_n(m) \eta \ra 
=  \lim_{n \to \infty} \la U(s_n) \xi, H(m) \eta \ra 
=  \la \xi, H(m) \eta \ra,\] 
and the boundedness of $H(M)\eta$ implies that  the convergence 
is uniform on $M$. This shows that the bounded function $H \: M \to B(\cH)$ 
is weakly holomorphic, hence holomorphic by \cite[Cor.~A.III.5]{Ne00}. 
Finally, the relation $U(f(m)) = U(g) U(h(m)) = U(g) H(m)$ 
implies that $U \circ f$ is holomorphic.   
\end{prf}

The following theorem is the main result of this section. 
It provides a mechanism to construct 
unitary endomorphisms of~$V$ 
by the inclusion $S^U_{\rm inv} \subeq S_V$. 
It implements the analytic continuation process from 
Theorem~\ref{thm:lo3.18} inside the Olshanski semigroup~$S^U$. 

\begin{thm} \mlabel{thm:4.17} {\rm(Inclusion Theorem)} 
Let $G$ be a $1$-connected Lie group with the involution $\tau_G$ 
and $\tau \in \Aut(\g)$ the induced automorphism. 
Further, let $(U, \cH)$ be a continuous antiunitary representation 
of $G \rtimes \{\id_G, \tau_G\}$ 
with discrete kernel and consider the standard subspace 
$V \subeq \cH$ specified by $J_V= U(\tau_G)$ and 
$\Delta_V = e^{2\pi i \partial U(h)}$ for some $h \in \g^\tau.$ 
Then 
\[ S^U_{\rm inv} \subeq S_V := \{ g \in G \: U(g) V \subeq V \}
\quad \mbox{ and } \quad (\g + i C_U)_{\rm inv} \subeq \L(S_V).\] 
\end{thm}

\begin{prf} We write $U \: S^U \to B(\cH), g\exp(ix) 
\mapsto U(g) e^{i\partial U(x)}$ for the canonical extension of the 
unitary representation $U$ to $S^U$ 
(Proposition~\ref{prop:ucprop-new}). 
For $s \in S^U_{\rm inv}$, 
we consider the bounded function 
\[ F \: \oline{\cS_\pi} \to B(\cH), \quad 
F(z) := U(\beta_z(s))\] 
which is defined because $\beta_z(s)\in S^U$ 
for $z \in \oline{\cS_\pi}$. We have 
\[ F(z+ t) = U(\exp th) F(z) U(\exp(-th)) 
= \Delta_V^{-it/2\pi} F(z) \Delta_V^{it/2\pi} \quad \mbox{ for } \quad 
t \in\R, z \in \oline{\cS_\pi},\] 
and $F$ is strongly continuous (Proposition~\ref{prop:ucprop-new}). 
That it is holomorphic on $\cS_\pi$ follows from 
Lemma~\ref{lem:autoholo} and 
the holomorphy of the map 
$\cS_\pi \to \Gamma_G(C_{U}), z \mapsto \beta_z(s)$ in 
the sense of Definition~\ref{def:holodisc}. 
We further note that $F(0) = U(s) \in U(G)$ is unitary and that 
\[ J_V F(0) J_V 
= J_V U(s) J_V 
=  U(\tau_G(s)) 
=  U(\beta_{\pi i}(s)) = F(\pi i).\] 
Now  Theorem~\ref{thm:lo3.18}(iv) 
implies that $U(s) = F(0) \in \cA_V$, 
and thus $s \in S_V$. 

The inclusion of the Lie wedges is an immediate consequence 
of $\L(S^U_{\rm inv}) = (\g + i C_U)_{\rm inv}$  
(Lemma~\ref{lem:wick-olsh}(iv)). 
\end{prf}

\begin{prob} Show that we actually have the equality 
$S^U_{\rm inv} = S_V$.  In the Germ Theorem (Theorem~\ref{thm:keylem3} 
below) 
we shall see that both subsemigroup do at least have the same germ, 
i.e., that there exists an $e$-neighborhood $\cU$ in $G$ with 
$S^U_{\rm inv} \cap \cU= S_V \cap \cU$. 
\end{prob}

\begin{rem} \mlabel{rem:3grad} 
(a) The construction in the preceding proof shows that, 
for $s\in S^U_{\rm inv}$,  the element 
$s^c := \beta_{\frac{\pi i}{2}}(s) = h \exp(ix)\in S^U$ satisfies 
\[ \hat U(s) = \alpha_{\frac{\pi i}{2}}(U(s)) = U(s^c) 
= U(h) e^{i \partial U(x)}.\] 
Therefore $\hat{U(s)}$ 
is injective with dense range (cf.\ Lemma~\ref{lem:regops}). 
If $S_V$ coincides with $S^U_{\rm inv}$, this implies that 
$\hat{S_V} \subeq B(\cH)_{\rm reg}$. 

(b) On the level of the Lie wedge $\L(S_V)$, we know that 
\[ (\g + i C_U)_{\rm inv} 
= (C_U \cap \fq_1) \oplus \fh_0 \oplus (-C_U \cap \fq_{-1}) \] 
is a cone with a rather simple structure and completely determined 
by the pair $(\tau, h)$ and the cone~$C_U$. 
In Section~\ref{sec:4}, we study $S_V$ from a different perspective 
and we shall see that this cone actually generates a $3$-graded Lie subalgebra
(Theorem~\ref{thm:main2} and Corollary~\ref{cor:3grad}). 
\end{rem}

The global structure of the semigroup 
$S^U_{\rm inv}$ is hard to analyze in the non-abelian 
context. In Subsection~\ref{subsec:5.1} we explain how to reduce 
the determination of this semigroup to the case 
where $e^{2\pi i \ad h} = \id_{\g_\C}$, i.e., $\ad h$ is diagonalizable 
with integral eigenvalues. 

\section{The 
subsemigroups $S_V$ in finite dimensional groups} 
\mlabel{sec:4}

As in Section~\ref{subsec:3.3}, we start in this section 
with an antiunitary representation $(U, \cH)$ of a 
finite dimensional graded Lie group $G \rtimes \{\id_G, \tau_G\}$, 
where $G$ is $1$-connected, 
and consider the standard 
subspace $V = V_\gamma$ determined by 
$J_V = U(\tau_G)$ and $\Delta_V^{-it/2\pi} = U(\exp th)$ 
with $h \in \fh$ as explained in the introduction. 
Under the assumption that $U$ has discrete kernel, we 
determine the Lie wedge of the closed subsemigroup 
\[ S_V = \{ g \in G \: U(g) V \subeq V \}\quad 
\mbox{ with the unit group} \quad G_V := \{ g \in G \: U(g) V =V \}.\] 
Our main result on $\L(S_V)$  (Theorem~\ref{thm:main2}) 
asserts that 
\[ \L(S_V) = (-C_U \cap \fq_{-1}(h)) \oplus \fh_0(h) \oplus 
(C_U \cap \fq_1(h))\] 
and that $\L(S_V)$ spans a $3$-graded Lie subalgebra of $\g$. 
This result is based on the Germ Theorem (Theorem~\ref{thm:keylem3}),  
asserting the existence of an  $e$-neighborhood $\cU$ in $G$ with 
$S^U_{\rm inv} \cap \cU= S_V \cap \cU$, which implies in particular 
that $\L(S_V) = \L(S^U_{\rm inv})$.

In Subsection~\ref{subsec:4.2} we discuss the unit group 
$G_V$ of $S_V$, and we discuss some examples in Subsection~\ref{subsec:4.3}.

\subsection{The Lie wedge of $S_V$} 
\mlabel{subsec:4.1}

The following theorem shows that the subsemigroups $S_V$ and $S^U_{\rm inv}$ 
of $G$ have the same germ, i.e., identical intersection with some 
$e$-neighborhood. 
\begin{thm} \mlabel{thm:keylem3} {\rm(The Germ Theorem)} 
If $\ker(U)$ is discrete, then there exists an $e$-neighborhood $\cU \subeq G$ 
such that $S_V \cap \cU = S_{\rm inv} \cap \cU$. 
\end{thm}

\begin{prf} In Theorem~\ref{thm:4.17} we have already seen that 
$S^U_{\rm inv} \subeq S_V$. Therefore it suffices to find $\cU$ such that 
$S_V \cap \cU \subeq S^U_{\rm inv}$. 
We write $\cH^\omega \subeq \cH$ for the subspace of analytic 
vectors of the representation $U$ of $G$ 
(cf.~\cite[Thm.~4]{Nel59}). 

\nin {\bf Step 1:} By \cite[Lemma~9.16]{HN93}, there exists a 
dense subspace $\cD \subeq \cH^\omega$ which is equianalytic in the sense 
that there exists an open 
convex circular $0$-neighborhood $\cW \subeq \g_\C$, 
such that the series 
\[ U^\eta(x) := \sum_{n = 0}^\infty \frac{1}{n!} (\dd U(x))^n \eta \] 
converges for $\eta \in \cD$ and $x \in \cW$ and 
defines a holomorphic function $U^\eta : \cW \to \cH$. 
It satisfies 
\begin{equation}
  \label{eq:tildeident}
U^\eta(x) = U(\exp x)\eta \quad \mbox{  for } \quad x \in \cW \cap \g 
\end{equation}
(see~\cite[\S 6]{Ne11} for more on analytic vectors).  

\nin {\bf Step 2:} We claim that 
\begin{equation}
  \label{eq:2l}
 U^\eta(ix) = e^{i \partial U(x)}\eta \quad \mbox{  for } \quad 
x \in \cW,
\end{equation}
where the right hand side has to be understood in terms of the 
measurable functional calculus for the selfadjoint operator $i\partial U(x)$. 
First we observe that the function 
\[ F \: \{z \in \C \: zx \in \cW \} \to \cH,\quad 
F(z) = U^\eta(zx) \] 
is holomorphic and coincides for $z \in [-1,1] \subeq \R$ with $U(\exp zx)\eta$. 
Let $\eps > 0$ be such that 
$\Omega := \{ a + i b \: |a| < \eps, -\eps < b < 1 + \eps\}$ satisfies 
$\Omega x\subeq \cW$. 
Then the continuous 
function $F(t) := U(\exp tx)\eta$ on $(-\eps,\eps)$ 
has an analytic continuation to $\Omega$, and this implies that 
it even extends to the strip $\Omega' := \R + (-\eps, 1 + \eps) i \subeq \C$, 
so that \cite[Lemma~A.2.5]{NO18} shows that 
$\eta \in \cD(e^{ti \partial U(x)})$ for $-\eps < t < 1 + \eps$. 
As the function $\Omega' \to \cH, 
z \mapsto e^{z \partial U(x)}$ is also holomorphic, we obtain 
$U^\eta(ix) = e^{i \partial U(x)}$ by analytic continuation. 

\nin {\bf Step 3:} Let $\cW' \subeq \cW \subeq \g_\C$ 
be an open convex $0$-neighborhood such that 
the Baker--Campbell Hausdorff product 
$x * y$ is defined by the convergent series for 
$x,y \in \cW'$ and defines a holomorphic function 
$\cW' \times \cW' \to \cW \subeq \g_\C$ (\cite[\S 9.2.5]{HN12}). 
We claim that 
\begin{equation}
  \label{eq:3l}
U^\eta(x_1 * x_2) 
= U(\exp x_1) U^\eta(x_2) \quad \mbox{ for } \quad 
x_1 \in \cW' \cap \g, x_2 \in \cW'
\end{equation}
if $\cW'$ is chosen small enough. 
As $U(\exp x_1)$ is unitary, both sides 
define $\cH$-valued functions holomorphic in $x_2$. 
Fixing $x_1 \in \g \cap \cW'$, both sides are 
equal by \cite[Lemma~6.7]{Ne11} if we choose $\cW'$ small enough. 
Here we use that the invariant subspace 
$U(\g_\C)\eta$ generated by $\eta$ under the derived representation 
is equianalytic (\cite[Prop.~6.6]{Ne11}). 

\nin {\bf Step 4:} For $x_1, x_2 \in \g \cap \cW'$ we 
now obtain with \eqref{eq:2l} and \eqref{eq:3l}
\begin{equation}
  \label{eq:4l}
U^\eta(x_1 * i x_2) = U(\exp x_1) e^{i \partial U(x_2)} \eta.
\end{equation}
Shrinking $\cW'$ if necessary, we may further assume that the map 
\[   (\cW' \cap \g) \times (\cW' \cap \g) \to \cW, \quad 
(x,y) \mapsto x * iy \] 
is a diffeomorphism onto an open subset $\cW''\subeq \cW$. 
As $\ker(U)$ is discrete, we may further assume that 
$U \circ \exp$ is injective on $\cW' \cap \g$. 

\nin {\bf Step 5:} Now let 
$\tilde \cW :=  \bigcap_{0 \leq y \leq \pi} e^{-yi \ad h} (\cW' \cap \cW'')$ 
and observe that, by the compactness of $[0,\pi]$,  
this is an open convex $0$-neighborhood in $\g_\C$.  
For $x \in \tilde\cW \cap \g$ with 
$\exp x \in S_V$, we then find an $\eps > 0$ such that 
\[  e^{z \ad h} x \in \cW' \cap \cW \quad \mbox{ for } \quad z \in \Omega'' 
:= \{ w = a + i b \in \C \: |a|< \eps, -\eps < b < 1 + \eps\}.\]
Then $\Omega'' \to \cH, z \mapsto U^\eta(e^{z \ad h}x)$ is a holomorphic function, 
and so is the continuous function 
\[ \oline{\cS_\pi} \to \cH, \quad z \mapsto \alpha_z(U(\exp x)) \eta\] 
on the open strip $\cS_\pi$ (Theorem~\ref{thm:lo3.18}). 
Since both functions coincide on the interval 
$(-\eps,\eps) \subeq \R$, we obtain by analytic continuation 
\begin{equation}
  \label{eq:star}
\alpha_{ti}(U(\exp x)) \eta = U^\eta(e^{ti \ad h}x)\quad \mbox{ for } \quad 
0 \leq t \leq \pi.
\end{equation}
As $e^{it \ad h}x \in \cW''$ for $0 \leq t \leq \pi$, 
we can write this element uniquely as 
$x_t * i y_t$ with $x_t, y_t \in \g \cap \cW'$. We now obtain with 
\eqref{eq:4l}
\[ \|U^\eta(x_t * i y_t)\| 
= \|U(\exp x_t) e^{i \partial U(y_t)}\eta\| 
= \|e^{i \partial U(y_t)}\eta\| \]
and 
\[ \|U^\eta(x_t * i y_t)\| = \|U^\eta(e^{it \ad h}x)\|
= \|\alpha_{it}(U(\exp x)) \eta\|
\leq \|\eta\| \] 
because $\|\alpha_{it}(U(\exp x))\| \leq \|U(\exp x)\| = 1$ 
by Theorem~\ref{thm:lo3.18}. Comparing both terms, we see that 
\[ \|e^{i \partial U(y_t)} \eta \| \leq \|\eta\| \quad \mbox{ for } \quad \eta \in \cD.\]  
As the operator $e^{i \partial U(y_t)}$ is selfadjoint, it is in particular closed. 
The above estimate shows that the closure of the restriction 
$e^{i \partial U(y_t)}\res_{\cD}$ is a bounded operator on $\oline\cD = \cH$. 
We conclude that $\|e^{i \partial U(y_t)}\| \leq 1$, and hence 
that $i\partial U(y_t) \leq 0$. This implies that $y_t \in C_U$ for 
$0 \leq t \leq \pi$. This in turn shows that 
\[\exp(e^{t i \ad h}x) = \exp(x_t * i y_t) = \exp(x_t) \exp(i y_t) \in S^U\] 
and therefore 
$\beta_{it}(\exp x) \in S^U$ exists for $t \in [0,\pi]$. 
To see that $\exp x \in S^U_{\rm inv}$, it remains to show that 
\begin{equation}
  \label{eq:taurel}
\beta_{\pi i}(\exp x) = \tau_G(\exp x) = \exp(\tau(x)).
\end{equation}
From Theorem~\ref{thm:lo3.18} we recall that 
\[ \alpha_{\pi i}(U(\exp x)) 
= J_V U(\exp x) J_V = U(\exp \tau(x)).\] 
We thus obtain for $t = \pi$ and $\eta \in \cD$: 
\[\|\eta\| =  \|\alpha_{\pi i}(U(\exp x))\eta\| 
= \|e^{i \partial U(y_\pi)}\eta\|,\] 
and since $i \partial U(y_\pi) \leq 0$, this leads to 
$\partial U(y_\pi)\eta = 0$. As $\ker(U)$ is discrete, it follows 
that $y_\pi = 0$, so that $e^{\pi i \ad h} x = x_\pi \in \g$. 
Now \eqref{eq:star} yields 
\[ U(\exp \tau(x))\eta = \alpha_{\pi i}(U(\exp x))\eta 
= U^\eta(e^{\pi i \ad h}x)
= U(\exp(e^{\pi i \ad h}x))\eta \quad \mbox{ for } \quad \eta \in \cD,\] 
which in turn leads to $U(\exp(e^{\pi i \ad h}x)) = U(\exp \tau(x))$. 
As $U \circ \exp$ is injective on $\cW' \cap \g$, 
we obtain $\beta_{\pi i}(\exp x) = \exp(e^{\pi i \ad h}x) = \tau_G(\exp x)$, 
and this finally proves that $\exp x \in S^U_{\rm inv}$ 
for $x \in \tilde \cW \cap \g$ with $\exp x \in S_V$. 
\end{prf}

The following corollary is a converse to Theorem~\ref{thm:4.17} 
on the level of infinitesimal generators. 
It follows immediately from the Germ Theorem (Theorem~\ref{thm:keylem3}), 
Lemma~\ref{lem:wick-olsh}(iv),  
and the observation that subsemigroups with identical germs 
have identical Lie wedges. 
 
\begin{cor} \mlabel{cor:keylem2} If $\ker(U)$ is discrete, then 
\[ \L(S_V) = (\g + i C_{U})_{\rm inv}
:= \{ x \in \g + i C_U \: 
e^{[0,\pi]i \ad h}x \subeq (\g + i C_U); e^{\pi i \ad h}x = \tau(x)\}.\] 
\end{cor}

Before we can prove the Structure Theorem, we need one more ingredient. 
We recall that a {\it standard pair} $(U,V)$ consists of a standard subspace 
$V \subeq \cH$ and a unitary one-parameter group $(U_t)_{t \in \R}$ 
satisfying $U_t V \subeq V$ for $t \geq 0$. 

\begin{prop} \mlabel{prop:trans-commute} 
Let $(G,\eps_G)$ be a finite dimensional graded Lie group 
and $(U,\cH)$ be an antiunitary representation of~$G$. 
Suppose that $(V,U^j)$, $j = 1,2$, are standard pairs for which 
there exists a graded homomorphism $\gamma \: \R^\times\to G$ 
and $x_1, x_2 \in \g$ such that 
\[ J_V = U(\gamma(-1)), \quad \Delta_V^{-it/2\pi} = U(\gamma(e^t)), 
\quad \mbox{ and }\quad 
U^j(t) = U(\exp t x_j), \quad t \in \R, j =1,2.\] 
Then the unitary one-parameter groups $U^1$ and $U^2$ commute. 
\end{prop}

In Subsection~\ref{subsec:parisstand} we describe an example 
showing that, without assuming that they come from a 
finite dimensional Lie group $G$, 
the two one-parameter groups $U^1$ and $U^2$ need not commute.

\begin{prf} The positive cone $C_U\subeq \g$ of the representation $U$ 
is a closed convex $\Ad(G)$-invariant cone. As we may w.l.o.g.\ assume that 
$U$ is injective, the cone $C_U$ is pointed. 

Writing $\Delta_V^{-it/2\pi} = U(\exp t h)$ and 
$U^j_t = U(\exp t x_j)$ with $h,x_1, x_2 \in \g$, we have 
$[h,x_j] = x_j$ for $j =1,2$ and $x_1, x_2 \in C_U$. 
If 
\[ \g_\lambda(h)= \ker(\ad h - \lambda \1) \]  
is the $\lambda$-eigenspace of $\ad h$ in $\g$, then 
$[\g_\lambda(h),\g_\mu(h)] \subeq \g_{\lambda + \mu}(h)$, so that 
$\g_+ := \sum_{\lambda >0} \g_\lambda(h)$ is a nilpotent Lie algebra. 
Therefore $\fn := (C_U \cap \g_+) - (C_U \cap \g_+)$ is a nilpotent 
Lie algebra generated by the pointed invariant cone $C_U \cap \g_+$.
By \cite[Ex.~VII.3.21]{Ne00}, $\fn$ is abelian. Finally 
$x_j \in C_U \cap \g_1(h) \subeq \fn$ implies that $[x_1,x_2] =0$.
\end{prf}

The following theorem not only provides an explicit description 
of the Lie wedge $\L(S_V)$, we also show that $\L(S_V)$ 
spans a $3$-graded Lie subalgebra $\g_{\rm red}$ of $\g$. 

\begin{thm} {\rm(Structure Theorem for $\L(S_V)$)} 
\mlabel{thm:main2}
If $\ker(U)$ is discrete, 
then 
\begin{equation}
  \label{eq:ls-form}
 \L(S_V) = C_- \oplus \fh_0(h) \oplus C_+, 
\quad \mbox{ where } \quad 
C_\pm =  \pm C_U \cap \fq_{\pm 1}(h).
\end{equation}
If $\fq_\pm := C_\pm - C_\pm$ are the linear subspaces generated by $C_{\pm}$, then 
$\L(S_V)$ spans 
the $3$-graded Lie subalgebra $\g_{\rm red} :=  \fq_- \oplus \fh_0(h) \oplus \fq_+$. 
\end{thm}

\begin{prf} From Corollary~\ref{cor:keylem2} 
we know that 
$(\g + i C_U)_{\rm inv} =  \L(S_V)$. 
Further, Corollary~\ref{cor:keylema-lie} implies that 
\[ (\g + i C_U)_{\rm inv} 
= (C_U \cap \fq_1(h)) \oplus \fh_0(h) \oplus (-C_U \cap \fq_{-1}(h)).\] 
This proves \eqref{eq:ls-form}. 
It follows in particular that $\fq_\pm = C_\pm - C_\pm \subeq \g_{\pm 1}(h)$. 
Proposition~\ref{prop:trans-commute} shows that the two 
subspaces $\fq_\pm$ of $\g$ are abelian. 
Further, $\fq_\pm \subeq \fg_{\pm 1}(h)$ implies that 
$[\fq_-,\fq_+] \subeq \fh \cap \g_0(h) = \fh_0(h)$,  
and  from Corollary~\ref{cor:sv1}  
we know that $\fh_0(h) = \L(S_V) \cap -\L(S_V)$.   
As the cone $\L(S_V)$ is a Lie wedge, the operators $e^{\ad x}$, $x \in \L(S_V) \cap \fh$ 
on $\fq$ preserve the cone $\L(S_V) \cap \fq$. This shows that 
$[[\fq_+, \fq_-], \fq_\pm] \subeq \fq_\pm,$ 
which implies that $\g_{\rm red}$ 
is a Lie subalgebra of $\g$. It is clearly $3$-graded 
by $\ad h$, and the restriction of $\tau$ to $\g_{\rm red}$ 
coincides with the restriction of $e^{\pi i \ad h}$. 
\end{prf}

\begin{cor} {\rm(The wedge $\L(S_V)$ in the $3$-graded case)} 
\mlabel{cor:3grad}
Suppose that $U$ has discrete kernel, and that 
\[ \g = \g_{-1}(h) \oplus \g_0(h) \oplus \g_1(h) 
\quad \mbox{ with }  \quad \tau = e^{\pi i \ad h},\] 
so that $\fh = \g_0(h)$ and $\fq = \g_1(h) \oplus \g_{-1}(h)$. 
Then 
\[ \L(S_V) = C_- \oplus \g_0(h) \oplus C_+, 
\quad \mbox{ where } \quad C_\pm = \g_{\pm 1}(h)\cap \pm C_U.\] 
\end{cor}

\subsubsection*{Independence of $S_V$ from $J_V$}

\begin{prop} {\rm(Independence of $S_V$ from $J_V$)} 
Let $(U^j,\cH)_{j=1,2}$ be antiunitary representations 
of the graded Lie group $(G,\eps_G)$ which coincide on $G_+$. 
Then, for every graded homomorphism $\gamma \: \R^\times \to G$, 
and the corresponding standard subspaces $V_\gamma^1$ and $V_\gamma^2$, we have 
\[ S_{V_\gamma^1} = S_{V_\gamma^2}.\] 
\end{prop}

\begin{prf} By \cite[Thm.~2.11]{NO17} the antiunitary representations 
$U^1$ and $U^2$ are equivalent because their restrictions to $G_+$ coincide. 
Hence there exists a unitary operator $\Phi \in \U(\cH)$ with \break 
$\Phi \circ U^1(g) = U^2(g) \circ \Phi$ for all $g \in G$. 
This implies in particular that $\Phi(V^1_\gamma) = V^2_\gamma$, 
and since $\Phi$ commutes with $U^1(G_+) = U^2(G_+)$, it follows that 
$S_{V^1_\gamma} = S_{V^2_\gamma}$. 
\end{prf}

\subsection{The unit group $G_V$}
\mlabel{subsec:4.2}

In the following we denote the centralizer of $x \in \g$ in $G$ 
by $C_G(x) := \{ g \in G \: \Ad(g)x = x\}$. 

\begin{lem} \mlabel{prop:1.42b} Suppose that $\ker(U)$ is discrete. 
The groups $G_V = S_V \cap C_{G_+}(h) \supeq C_{G_+}(h)^{\tau_G}$ 
have  the same Lie algebra 
$\g_V = \L(S_V) \cap \g_0(h) = \fh_0(h)$. They coincide 
if $U$ is injective. 
\end{lem}

\begin{prf} If $g \in G_+$ satisfies $\Ad(g) h = h$, i.e., $g \in C_{G_+}(h)$, 
then the unitary operator 
$U(g)$ commutes with $\Delta_V$. Therefore  Corollary~\ref{cor:degree0}  
implies that $U(g) V \subeq~V$ is equivalent to 
$U(g) V =~V$. If 
$g \in C_{G_+}(h)^{\tau_G}$, then 
$U(g)$ commutes with $J_V$ and $\Delta_V$, so that $U(g)V = V$ 
(Lemma~\ref{lem:1.1}). 
This shows that 
\[ C_{G_+}(h)^{\tau_G} \subeq G_V  = S_V \cap C_{G_+}(h).\] 
If, in addition, $U$ is injective, then $U(g) \in G_V$ implies that 
$U(g)$ commutes with $U(\tau_G) = J_V$, and therefore $\tau_G(g) = g$. 

For the Lie algebras of these groups, we obtain 
\[  \g_V = \L(S_V) \cap -\L(S_V) = \L(S_V) \cap \g_0(h) \supeq \fh_0(h).\] 
Since the kernel of $U$ is discrete and the derived 
representation is injective, the fact that every 
$x \in \g_V$ generates a unitary one-parameter 
group commuting with $J_V = U(\tau_G)$ (Lemma~\ref{lem:1.1}) 
implies that $\tau(x) = x$, i.e., $x \in \fh$. We conclude that 
$\g_V = \fh_0(h)$. This proves the assertion on the Lie algebras. 
\end{prf}

\begin{ex}\mlabel{ex:sl2}  (Inequality in Lemma~\ref{prop:1.42b}) 
We consider the group 
$G_+ = \tilde\SL_2(\R)$ whose center is $Z(G_+) \cong \Z \cong \pi_1(\PSL_2(\R))$. 
Here the fundamental group of $\PSL_2(\R)$ it is generated by the loop 
obtained from the inclusion 
$\PSO_2(\R) \into \PSL_2(\R)$. 
Let $\tau_G \in \Aut(G_+)$ be the involution given on the Lie algebra level 
by $\tau\pmat{a & b \\ c & d} =   \pmat{a & -b \\ -c & d}$, and observe
that it induces the map $\tau_G(z) = z^{-1}$ on $Z(G_+)$. 

Now consider an antiunitary representation $(U,\cH)$ of $\PGL_2(\R)$, 
so that the corresponding representation of $G := G_+ \rtimes \{\id,\tau_G\}$ 
has kernel $Z(G_+)$. We therefore have $Z(G_+) = \ker(U) \subeq G_V$ 
for every $V \in \Stand(\cH)$. On the other hand, 
$Z(G_+) \subeq C_{G_+}(h)$, but $Z(G_+)$ 
it is not pointwise fixed by  $\tau_G$. We therefore 
have a proper inclusion $C_{G_+}(h)^{\tau_G} \into G_V$ in 
Lemma~\ref{prop:1.42b}. 
\end{ex}

\begin{rem} (Degenerate cases) 
\mlabel{rem:4.7} (a) If $\fq = \{0\}$ and $G_+$ is connected, then 
$\tau_G = \id_G$, so that Lemma~\ref{prop:1.42b} and 
Corollary~\ref{cor:sv1} imply that 
$S_V = G_V = C_{G_+}(h)$. 

\nin (b) If $\L(S_V) \cap \fq = \{0\}$, then the Structure Theorem~\ref{thm:main2} 
 implies that $\L(S_V) = \fh_0(h)$ is a Lie subalgebra of $\fh$. 
Now  \cite[Thm.~IV.2.11]{HHL89} implies that the group 
$G_V$ is ``isolated'' in $S_V$, i.e., there exists an $e$-neighborhood 
$\cU \subeq G$ with $\cU \cap S_V \subeq G_V$. Here we may w.l.o.g.\ 
assume that $\cU = \cU G_V$ is a tubular neighborhood of $G_V$. 
\end{rem}

\subsection{Some examples} 
\mlabel{subsec:4.3} 

\begin{ex} (a) (The affine group) For the graded group 
$G := \Aff(\R)= \R \rtimes \R^\times$ and representations of 
the form $U(b,a) := e^{ibP}  U^V(a)$, $V \in \Stand(\cH)$, 
we can use Theorem~\ref{thm:main2} to calculate the semigroup 
$S_V$. Since $S_V$ contains all positive dilations $(0,a)$, $a > 0$, 
this closed subsemigroup of $\R \rtimes \R^\times_+$ 
is of the form $S_V = C \rtimes \R^\times_+$, where 
$C = (\R \times \{1\}) \cap S_V$ is a closed additive subsemigroup 
of $\R$ invariant under multiplication with positive scalars. 
This leaves only the possibilities $C = \{0\}, [0,\infty)$ or $(-\infty,0]$.
Comparing with the Structure Theorem~\ref{thm:main2}, we obtain 
\begin{equation}
  \label{eq:affR-case}
  S_V = C \rtimes \R^\times, \quad \mbox{ where } \quad 
C := \{ x \in \R \:  x P \geq 0\}.
\end{equation} 

This case can also be derived from the standard 
subspace version of the Borchers--Wiesbrock Theorem 
(\cite[\S 3.2]{Lo08} and \cite[Thms.~3.13, 3.15]{NO17}). 

(b) (Higher dimensional dilation groups) More generally, 
we consider a group of the form $G = E \rtimes_\alpha \R^\times$, 
where the homomorphism $\alpha \: \R^\times  \to \GL(E)$ satisfies 
$\alpha(r) = r \1$ for $r > 0$. Then $\tau_E := \alpha(-1)$ is an involution 
and we write $E = E^+ \oplus E^-$ for the corresponding eigenspace 
decomposition. 

Let $(U,\cH)$ be an antiunitary representation 
of $G$, and  consider the standard subspace 
$V \in \Stand(\cH)$ with $U^V(r) = U(0,r)$ for $r \in \R^\times$. 
Then we also have 
\[ S_V = C \rtimes \R^\times_+, \] 
where $C \subeq E$ is a closed subsemigroup containing $0$ which is 
invariant under multiplication with positive scalars, hence a closed convex cone. 
As $\fh_0(h) = \{0\} \times \R$ and $\g_1(h) = E$, our 
Structure Theorem implies that 
$C = C_U \cap E^-.$ 
On the other hand, \cite[Thm.~3.15]{Lo08} implies that 
$E \cap C_U \subeq E^-$, so that
\[ C = C_U \subeq~E^-.\] 
Note that we cannot apply (a) directly to the one-dimensional subspaces 
of $E$ because we did not assume that $\alpha(-1) = -\id_E$. 
\end{ex}

\begin{ex} (More general $\R^\times$-actions) We consider a group of the form 
$G = E \rtimes_{\alpha} \R^\times$, so that 
$\tau_E := \alpha(-1)$ is an involutive automorphism of $E$. Accordingly, 
we write $E = E^+ \oplus E^-$ with $E^\pm = \ker(\tau_E \mp \1)$ for the 
$\tau_E$-eigenspace decomposition. We then have 
$\g = E \rtimes \R h$ with $\fq = E^-$ and $\fh = E^+ \rtimes \R h$. 
As $S_V$ contains $\{0\} \rtimes \R^\times$, we have 
\begin{equation}
  \label{eq:4.10}
S_V = (S_V \cap E) \rtimes_\alpha \R^\times_+, 
\end{equation}
where $S_V \cap E$ is a closed subsemigroup of $E$, invariant under 
$\alpha(\R^\times_+)$. 
We know from Theorem~\ref{thm:main2} that 
\[ \L(S_V) = (E^-_1(h) \cap C_U) \oplus (-E^-_{-1}(h) \cap C_U) \oplus 
(E^+_0(h) \oplus \R h),\] 
where 
\[ \L(S_V) \cap E^- = (E^-_1(h) \cap C_U) \oplus (-E^-_{-1}(h) \cap C_U) \] 
is a pointed convex cone determined by the positive cone $C_U$ of~$U$, and 
\[ \L(S_V) \cap - \L(S_V) = E^+_{0}(h) \oplus \R h\] 
is a Lie subalgebra. Here we can even use Subsection~\ref{subsec:3.1} 
to determine the subsemigroup $S^U_{\rm inv}$ of $S_V$. 
From $S^U_{\rm inv} = (S^U_{\rm inv} \cap E) \rtimes \R^\times_+$ 
and 
\[ S^U_{\rm inv} \cap E 
= (E + i (C_U \cap E))_{\rm inv} = \L(S_V) \cap E,\] 
we obtain 
\begin{equation}
  \label{eq:dag}
S^U_{\rm inv} = (\L(S_V) \cap E) \rtimes \R^\times_+,
\end{equation}
so that $S^U_{\rm inv}$ is the maximal infinitesimal generated subsemigroup 
of $S_V$. Presently we do not know if we always have $S^U_{\rm inv} = S_V$ 
for this class of groups, but this is work in progress 
(\cite{Ne19}).  
\end{ex}

\begin{ex}
Suppose that $\g$ is a simple real Lie algebra, 
that $G = G_+ \rtimes \{\id,\tau_G\}$, where the group $G_+$ is connected, 
and that $(U,\cH)$ is an antiunitary representation with 
non-zero positive cone~$C_U$ and discrete kernel. 
This already implies that 
$\g$ is quite special, it has to be a hermitian Lie algebra 
(see \cite{Ne00} for details and a classification). 
As $J := U(\tau_G)$ is antiunitary, 
$-\tau(C_U) = C_U$ by Remark~\ref{rem:4.15}. 
We pick $h \in \fh = \g^\tau$ and consider the corresponding 
semigroup~$S_V$. 

If $\fq = \{0\}$, then $\tau = \id_\g$, so that 
$U(G_+) \subeq \U(\cH)^J$ implies that $S_V$ is a group, 
namely the centralizer of $\Delta_V$ in $G_+$ (Remark~\ref{rem:4.7}(a)). 
We may therefore exclude this case and assume that $\tau \not=\id_\g$. 

The Structure Theorem (Theorem~\ref{thm:main2}) shows that 
\[ \L(S_V) = C_- \oplus \fh_0(h) \oplus C_+, 
\quad \mbox{ where } \quad 
C_\pm =  \pm C_U \cap \fq_{\pm 1}(h).\] 

In general, this cone may be rather small, but we know from 
Theorem~\ref{thm:main2} that it spans a $3$-graded Lie algebra $\g_{\rm red}$. 
If $\g = \g_{\rm red}$, 
then $\g$ itself is $3$-graded, hence a hermitian Lie algebra 
of tube type, i.e., the conformal Lie algebra of a euclidean 
Jordan algebra (see \cite[\S 3]{Ne18} for more details). 
In this case and for the centerless group $G$ with Lie algebra $\g$, 
we have determined the semigroup $S_V$ in 
\cite[Thms.~3.8, 3.13]{Ne18}: It coincides with the product set 
\[ S_V = \exp(C_+) (G_+^{\tau_G}) \exp(C_-).\] 
\end{ex}

\section{Perspectives}
\mlabel{sec:5} 

\subsection{Relations to von Neumann algebras} 
\mlabel{subsec:5.0}

We already mentioned in the introduction that the interest 
in the semigroups $S_V$ of endomorphisms of standard subspaces 
stems to some extent from their correspondence  to endomorphisms  
of von Neumann algebras 
in the context of the theory of local observables 
(\cite{Ha96}). We now provide some more details on these applications. 

We recall the notion of a 
{\it Haag--Kastler net} of $C^*$-sub\-algebras $\cA(\cO)$ 
of a $C^*$-algebra $\cA$,  
associated to regions $\cO$ in $d$-dimensional Minkowski space $\R^{1,d-1}$. 
The algebra $\cA(\cO)$ is interpreted as observables  
that can be measured in the ``laboratory'' $\cO$. 
Accordingly, one requires {\it isotony}, i.e., that $\cO_1 \subeq \cO_2$ implies 
$\cA(\cO_1) \subeq \cA(\cO_2)$ and that the $\cA(\cO)$ generate $\cA$. 
Causality enters by the {\it locality} assumption that 
$\cA(\cO_1)$ and $\cA(\cO_2)$ commute if 
$\cO_1$ and $\cO_2$ are space-like separated, i.e., cannot correspond 
with each other. 
Finally one assumes an action 
$\sigma \: P(d)_+^{\up} \to \Aut(\cA)$ 
of  the connected Poincar\'e group such that 
$\sigma_g(\cA(\cO)) = \cA(g\cO)$. Every Poincar\'e invariant state $\omega$ 
of the algebra $\cA$ now leads by the GNS construction to a covariant 
representation $(\pi_\omega, \cH_\omega, \Omega)$ of $\cA$, and hence 
to a net $\cM(\cO) := \pi_\omega(\cA(\cO))''$ of von Neumann algebras 
on $\cH_\omega$. Whenever $\Omega$ is cyclic and separating for 
$\cM(\cO)$,  we obtain modular objects 
$(\Delta_\cO, J_\cO)$. This connection between the 
Araki--Haag--Kastler theory of local observables  
and modular theory leads naturally to antiunitary group representations 
(cf.~\cite[\S 5]{NO17} and the introduction). 

Let, more generally, $(G,\eps_G)$ be a graded group of spacetime 
symmetries, where $\eps_G(g) = 1$ means that $g$ preserves time orientation 
and $\eps_G(g) = -1$ that it reverses time orientation; 
a typical example is the Poincar\'e group $P(d)$ with 
$P(d)_+ = P^\uparrow(d)$. Then covariant representations of 
Haag--Kastler nets lead to families $\cM(\cO)$ of von Neumann algebras 
and antiunitary representations $U \: G \to \AU(\cH)$ satisfying 
\[ U(g) \cM(\cO) U(g)^{-1} = \cM(g\cO).\] 
If the vacuum vector $\Omega\in \cH$ is fixed by  $U(G)$ and 
$\Omega$ is cyclic and separating for the von Neumann algebra 
$\cM(\cO)$, and $U(G)$ contains the corresponding modular 
conjugation $J$ and the one-parameter group $(\Delta_\cO^{it})_{t \in \R}$, 
then we are in the situation mentioned in the introduction, 
and we obtain information on  the subsemigroup 
\[ S_{V_\cO} \supeq 
S_\cO := \{ g \in G_+ \: U(g) \cM(\cO) U(g)^{-1} \subeq \cM(\cO)\}.\] 

Theorem~\ref{thm:main2} implies that the Lie wedge $\L(S_V)$ spans a $3$-graded 
Lie subalgebra $\g_{\rm red}$ such that the corresponding 
$3$-graded subgroup $G_{\rm red} \subeq G$ has the property 
that $S_V \cap G_{\rm red}$ has interior points and that the modular 
conjugation and the modular group also come from $U(G_{\rm red})$. 

\begin{ex}
  \begin{footnote}
{We thank Yoh Tanimoto for the discussion that led to this example.}
  \end{footnote}
(a) It is important to observe that, in the situation 
described in the introduction, where $\Omega$ is a cyclic separating 
unit vector for the von Neumann algebra $\cM$ and 
\[ V = \oline{\{ M\Omega \: M = M^* \in \cM\}},\] 
the inclusion 
\[ S_\cM 
= \{ g \in G \: g \cM g^{-1} \subeq \cM\} 
\subeq S_V = \{ g \in G \: g V\subeq V \} \] 
may be proper. 

To see an example, we consider the Hilbert space 
$\cH := B_2(\C^n)$ of Hilbert--Schmidt operators on $\C^n$ 
with the scalar product $\la A, B \ra := \tr(A^*B)$. 
By matrix multiplications from the left, we obtain a von Neumann 
subalgebra $\cM \subeq B(\cH)$, isomorphic to $M_n(\C)$, and its 
commutant $\cM'$ consists of right multiplications. 
The unit vector $\Omega := \frac{1}{\sqrt n}\1_n$ is cyclic and separating, and 
the corresponding standard subspaces for $\cM$ and $\cM'$ coincide with the space 
\[ V_\cM = V_{\cM'} = \Herm_n(\C)\] 
of hermitian matrices. Now $\theta(A) := A^\top$ defines a unitary 
operator on $\cH$ preserving $V_\cM = V_{\cM'}$ and 
satisfying $\theta \cM \theta^{-1} = \cM'$. 
For $G = \U(\cH)$, we therefore obtain $S_V \not= S_\cM$. 

(b) In the situation above, when $\cM$ is given, 
the $G$-orbit of $\cM$ in the space of von Neumann subalgebras of $B(\cH)$ 
can be identified with the homogeneous space 
$G/G_\cM$, and similarly, $G/G_V \into \Stand(\cH), g G_V \mapsto gV$ is an 
embedding. The discrepancy between both spaces comes from the 
potential non-triviality of the action of the stabilizer group $G_V$ 
on the von Neumann algebra~$\cM$. 

Related questions have been analyzed by Y.~Tanimoto in \cite{Ta10}. 
He refines the picture by considering the closed convex cone 
\[ V_\cM^+ = \oline{\{ M\Omega \: 0 \leq M = M^* \in \cM\}} \subeq V_\cM,\] 
which leads to the inclusions 
\[ S_\cM \into S_{V_\cM^+} = \{ g \in G \: g V_\cM^+ \subeq V_{\cM}^+\} 
\subeq S_{V_\cM}.\] 
Here the semigroup $S_{V_{\cM}^+}$ appears to be much closer to $S_\cM$ 
than $S_V$; see in particular \cite[Thm.~2.10]{Ta10}. 
In this context it is also interesting to note that the map 
\[ V_\cM^+ \to \cM_*^+, \quad \xi \mapsto \omega_\xi, \qquad 
\omega_\xi(M) = \la \xi, M\xi \ra \] 
is a homeomorphism by \cite[Thm.~1.2]{Ko80}. Accordingly, 
every element $g \in S_{V_\cM^+}$ induces a continuous map on $\cM_*^+$. 
\end{ex}

\begin{ex} In many situations arising in QFT, the 
group $G$ is the Poincar\'e group $P(d) \cong \R^{1,d-1} \rtimes \OO_{1,d-1}(\R)$ 
acting by affine isometries on $d$-dimensional Minkowski space $\R^{1,d-1}$. 
We define a grading on $P(d)$ by time reversal, i.e., 
$\eps_G(v,g) = \eps(g)$ and $g(V_+) = \eps(g) V_+$ for the upper open 
light cone $V_+ := \{ (x_0, \bx) \in \R^{1,d-1} \: x_0 > 0, 
x_0^2 > \bx^2\}$. 

The generator $h \in \so_{1,d-1}(\R)$ of the Lorentz boost on the 
$(x_0,x_1)$-plane 
\[ h(x_0,x_1, \ldots, x_{d-1}) = (x_1, x_0, 0,\ldots, 0).\] 
It satisfies $e^{2\pi i \ad h} = \1$, and 
$\tau := e^{\pi i \ad h}$ defines an involution on the 
Poincar\'e--Lie algebra $\fp(d)$, acting on $\R^{1,d-1}$ by 
\[ \tau(x_0, x_1, \ldots, x_{d-1}) = (-x_0, -x_1, x_2, \ldots, x_{d-1})\] 
on $\R^{1,d-1}$. 

For any positive energy representation of $P(d)$ with discrete kernel, we 
then have $C_U = \oline{V_+}$, because 
this is, up to sign, the only non-zero pointed
invariant cone in the Lie algebra $\fp(d)$ (for $d > 2$). 
Therefore the Lie wedge of the corresponding semigroup 
$S_V$ associated to the standard subspace determined by 
the triple $(U,\tau,h)$ is given by 
\[\L(S_V) = \fh_0(h) \oplus (\fq_1(h) \cap C_U) \oplus (\fq_{-1}(h) \cap -C_U)\] 
(Theorem~\ref{thm:main2}).
Here $\fh_0(h) = \g_0(h)$ is the centralizer of the Lorentz boost: 
\[ \g_0(h)=  (\{(0,0)\} \times \R^{d-2}) \rtimes (\so_{1,1}(\R) \oplus \so_{d-2}(\R)) 
\cong (\R^{d-2} \rtimes \so_{d-2}(\R)) \oplus \R h, \] 
and, for $\fq_j := \fq_j(h)$: 
\[ \fq_1 \cap C_U = \R (e_0 + e_1) \cap \oline{V_+} = \R_+ (e_1 + e_0) 
\quad \mbox{ and } \quad 
\fq_{-1} \cap (-C_U) = \R (e_0 - e_1) \cap -\oline{V_+} = \R_+ (e_1-e_0).\] 
Therefore $\L(S_V)$ coincides with the Lie wedge of the semigroup 
\[ S_{W_R} := \{ g\in P(d)_+ \:  gW_R \subeq W_R \},\] 
where $W_R := \{ x \in \R^{1,d-1} \: x_1 > |x_0|\}$ is the open 
right wedge (see also \cite[Lemma~4.12]{NO17}).  
\end{ex}

The starting point for the development that led to 
fruitful applications of modular theory in QFT 
was the Bisognano--Wichmann Theorem, asserting that  
the modular automorphisms $\alpha_t(M) = \Delta^{-it/2\pi}M \Delta^{it/2\pi}$ 
associated to the algebra $\cM(W_R)$ of observables corresponding 
to the right wedge $W_R$ in Minkowski space 
are implemented by the unitary action of a 
one-parameter group of Lorentz boosts preserving $W_R$. 
This geometric implementation of modular automorphisms in terms of 
Poincar\'e transformations was an important first step in a 
rich development based on the work of Borchers and Wiesbrock 
in the 1990s \cite{Bo92, Bo95, Bo97, Wi92, Wi93, Wi93c}.
They managed to distill the abstract essence from the Bisognano--Wichmann 
Theorem which led to a better understanding of the 
basic configurations of von Neumann algebras 
in terms of half-sided modular inclusions 
and modular intersections. 
In his survey \cite{Bo00}, Borchers described how 
these concepts have revolutionized quantum field theory. 
Subsequent developments can be found in \cite{Ar99, BGL02, Lo08, LW11, 
JM18, Mo17}.

\subsection{Pairs of standard pairs} 
\mlabel{subsec:parisstand}

For $V \in \Stand(\cH)$, one may expect that one-parameter 
groups $U^1$ and $U^2$, for which $(V,U^j)$ form a standard pair, 
commute. By Proposition~\ref{prop:trans-commute} this is true 
if they both come from an antiunitary representation 
of a finite dimensional Lie group. 
The following example shows that this is not true in 
general, not even if the two one-parameter groups are conjugate under the 
stabilizer group~$\U(\cH)_V$. 

\begin{ex} On $L^2(\R)$ we consider the selfadjoint operators 
\[ (Qf)(x) = xf(x) \quad \mbox{ and }\quad (Pf)(x) = i f'(x),\]
satisfying the canonical commutation relations $[P,Q] = i \1$. 
For both operators, the Schwartz space $\cS(\R)$ is a core. 
Actually it is the space of smooth vectors for the representation 
of the $3$-dimensional Heisenberg group generated by the corresponding 
unitary one-parameter groups 
\[ (e^{itQ}f)(x) = e^{itx}f(x) \quad \mbox{ and }\quad (e^{itP}f)(x) = f(x-t).\]
Since $e^{ix^3}$ is a smooth function for which all derivatives grow at most 
polynomially, it defines a continuous linear operator on $\cS(\R)$ 
(\cite[Thm.~25.5]{Tr67}). Therefore the unitary operator $T := e^{iQ^3}$ 
maps $\cS(\R)$ continuous onto itself, and 
\[ \tilde P := TPT^* = e^{iQ^3} P e^{-iQ^3} \] 
is a selfadjoint operator for which $\cS(\R)$ is a core. 
For $f \in \cS(\R)$, we obtain 
\[ (\tilde P f)(x) 
= i e^{ix^3} \frac{d}{dx} e^{-ix^3}f(x) 
= i (-i 3x^2 f(x) + f'(x)),\] 
so that $\tilde P = P + 3 Q^2.$

The two selfadjoint operators $Q$ and $e^P$ are the infinitesimal generators 
of the irreducible antiunitary representation of $\Aff(\R) = \R \rtimes \R^\times$, 
given by 
\[ U(b,e^t) = e^{ib e^P} e^{it Q} \quad \mbox{ and }\quad 
(U(0,-1) f)(x) = \oline{f(-x)}.\] 
Accordingly, the pair $(\Delta, J)$ with 
\[ \Delta = e^{-2\pi Q} \quad \mbox { and }\quad J = U(0,-1) \] 
specifies a standard subspace $V$ which combines with 
$U^1_t := e^{it e^P}$ to an irreducible standard pair~$(V,U^1)$. 
The unitary operator $T$ commutes with $\Delta$ and with $J$ because 
$JQJ = -Q$, so that $T(V) = V$. Therefore the unitary one-parameter group 
$U^2_t :=  e^{iQ^3} U^1_t e^{-iQ^3}=  e^{it e^{\tilde P}}$ 
also defines a standard pair $(V, U^2)$. These two one-parameter groups 
do not commute because otherwise the selfadjoint operators 
$P$ and $P + 3 Q^2$ would commute in the strong sense, hence in particular 
on their core $\cS(\R)$. 
\end{ex}

\subsection{Classification problems} 
\mlabel{subsec:5.1a}

In the light of our results on the structure of the Lie wedges 
$\L(S_V)$, one would like to classify all situations, 
where these cones generate the Lie algebra $\g$. 
As this requires $\g$ to be $3$-graded by $\ad h$ with 
$\tau = e^{\pi i \ad h}$, we have to consider 
Lie groups $G_+$ with Lie algebra $\g$ 
and $\ad$-diagonalizable elements $h \in \g$ with 
$\Spec(\ad h) \subeq \{-1,0,1\}$. 
Then we have to study  unitary representations of $G_+$ 
extending to antiunitary representations of $G = G_+ \rtimes \{\id_G, \tau_G\}$ 
in such a way that the two cones $\fg_{\pm 1}(h) \cap C_U$ generate $\g_{\pm 1}(h)$. 
Then the ideal 
$\g_1 := [\fq,\fq] \oplus \fq$ generated by $\fq$ is  contained 
in $C_U - C_U$, so that the cone $C_{U^1}$ for the restriction 
$U^1 := U\res_{G^1}$ is generating. 

Since we expect the semigroup $S_V$ to be adapted to 
any direct integral decomposition into irreducible representations, 
the main point is to understand the irreducible 
representations. 
For the normal subgroup 
$G^1$ we thus have to study irreducible antiunitary representations $U^1$ 
for which the cone $C_{U^1}$ is pointed and generating. 
Up to the extendability question from $G^1_+$ to $G^1$ 
(cf.~\cite[Thm.~2.11(c),(d)]{NO17}), we are then dealing with 
unitary highest weight modules, whose classification theory can be found in 
\cite[\S X.4]{Ne00}. So the first steps in a classification 
should start with a faithful unitary highest weight representation 
$(U_\lambda,\cH)$ of a one-connected Lie group $G^1_+$ and a 
derivation $D \in \der(\g_1)$ satisfying $D^3 = D$, such that 
$C_U \cap \fg_{\pm 1}(D)$ generates $\fg_{\pm 1}(D)$. 
Then $\g = \g_1 \rtimes_D \R$ is a Lie algebra to which our results apply.

\subsection{Global information on the semigroup $S_V$}
\mlabel{subsec:5.1}

We recall the context of Theorem~\ref{thm:4.17} 
with the semigroup $S^U = G \exp(i C_U)$ 
on which the analytic extension of the unitary representation 
$(U,\cH)$ of the $1$-connected Lie group $G$ lives, and the subsemigroup 
$S^U_{\rm inv} \subeq S_V$ which has the same germ as $S_V$ (Theorem~\ref{thm:keylem3}). 
Therefore the picture is very clear for the Lie wedges, but 
the global semigroup $S_V$ and $S^U_{\rm inv}$ may be more complicated and not 
even generated by their one-parameter subsemigroups. 
It would be interesting to understand the structure 
of the subsemigroup $S^U_{\rm inv} \subeq G$ better, 
but this problem is quite intricate as well. 
However, below we shall see that it reduces to the situation 
where $e^{2\pi i \ad h}= \1$, which is a 
non-abelian analog of Lemma~\ref{lem:4.13a}.

\begin{lem} \mlabel{lem:zetacomplex}
 Consider the $1$-connected complex Lie group $G_\C$ with Lie algebra $\g_\C$ 
and the two connected Lie subgroups 
$G := G_\C^\sigma$ and $G^c := G_\C^{\sigma^c}$, 
where $\sigma$ and $\sigma^c$ are the two antiholomorphic 
involutions of $G_\C$ for which the derivative in $e$ is complex conjugation 
with respect to $\g$ and $\g^c$, respectively. 
Then the following assertions hold: 
\begin{itemize}
\item[\rm(i)] $\tau_{G_\C} = \sigma\sigma^c = \sigma^c \sigma$ 
is the holomorphic involution integrating the complex linear 
extension of $\tau$ to $\g_\C$. 
\item[\rm(ii)] For $\zeta_G := \beta_{\pi i/2} \in \Aut(G_\C)$ and 
$g \in G$, we have $\zeta_G(g) \in G^c$ if and only if 
$\beta_{\pi i}(g) = \tau_G(g)$, and this implies that 
$\zeta_G^4(g) = g$, so that $\zeta_G^{-1}(G^c) \cap G \subeq \Fix(\zeta_G^4)$. 
\item[\rm(iii)] For elements of the form 
$g^c = h \exp(x) \in G^c$ with $h\in H^c := (G^c)^{\tau_G}$ and $x \in i \fq$ 
with $\Spec(\ad x) \subeq \R$, 
we have $g^c\in \zeta_G(G)$ if and only if 
$h \in \zeta_G(G)$ and $x \in \zeta(\g)$. 
If this is the case, then $e^{\pi i \ad h}x = -x$ and $\beta_{\pi i}(h) = h$. 
\end{itemize}
\end{lem}

\begin{prf} (i) follows by inspection of the differentials. 

\nin (ii) For the automorphisms $\beta_z \in \Aut(G_\C)$ 
with differential $e^{z \ad h}$, we have 
\[ \sigma \circ \beta_z = \beta_{\oline z} \circ \sigma \quad \mbox{ and } \quad 
\sigma^c \circ \beta_z = \beta_{\oline z} \circ \sigma^c \quad \mbox{ for } \quad 
z \in \C.\]
For $z = \pi i/2$, we obtain in particular 
\[ \sigma \circ \zeta_G = \zeta_G^{-1} \circ \sigma \quad \mbox{ and } \quad 
 \sigma^c \circ \zeta_G = \zeta_G^{-1} \circ \sigma^c.\] 
Now let $g \in G$. The condition $\zeta_G(g) \in G^c$ is by 
$\sigma^c(g) = \sigma^c \sigma(g) = \tau_G(g)$ equivalent to 
\[ \zeta_G(g) 
{\buildrel {!} \over =}  \sigma^c(\zeta_G(g)) = \zeta_G^{-1}(\sigma^c(g)) 
= \zeta_G^{-1}(\tau_G(g)), \]
hence to $\tau_G(g) = \zeta_G^2(g) = \beta_{\pi i}(g).$
If this condition is satisfied, then 
\[ g = \tau_G(\tau_G(g)) = \tau_G(\zeta_G^2(g)) = \zeta_G^2(\tau_G(g)) 
= \zeta_G^4(g).\] 

\nin (iii) If $h \in \zeta_G(G)$ and $x \in \zeta(\fg)$, then we clearly have 
$h \exp x \in \zeta_G(G)$. 

Suppose, conversely, 
that $g^c = h \exp x$ with $h \in H^c$ and $x \in i\fq$ with 
$\Spec(\ad x) \subeq \R$ satisfies $g \in \zeta_G(G)$. 
As $\zeta_G$ commutes with $\tau_G$ and the group 
$G$ is invariant under $\tau_G$, the group 
$\zeta_G(G)$ is also $\tau_G$-invariant. Hence 
$g^c \in \zeta_G(G)$ implies $g^\sharp \in \zeta(G)$ and thus also 
$g^\sharp g = \exp 2x \in \zeta_G(G)$. 
The latter condition can be written as 
$\exp(2\sigma\zeta^{-1}(x)) = \exp(2 \zeta^{-1}(x)).$ 
Since $\ad x$ has real spectrum and $G_\C$ is simply connected, 
we obtain with Lemma~\ref{lem:expinj} that 
$\sigma \zeta^{-1}(x) = \zeta^{-1}(x)$, i.e., $x\in \zeta(\g)$.
This in turn implies that $h \in \zeta_G(G)$. 

From (ii) we now obtain $\zeta(x) = \zeta^2(\zeta^{-1}(x)) 
= \tau(\zeta^{-1}(x)) = \zeta^{-1}(\tau(x)) = - \zeta^{-1}(x)$, 
hence $\zeta^2(x) = -x$. We likewise get 
$\zeta_G(h) = \zeta_G^2(\zeta_G^{-1}(h)) = \tau_G(\zeta_G^{-1}(h)) 
= \zeta_G^{-1}(\tau_G(h)) = \zeta_G^{-1}(h)$, and therefore 
$\zeta_G^2(h) = h$. 
\end{prf}

The following proposition 
reduces the determination of 
$S^U_{\rm inv}$ to the case where $\zeta^4 = \1$, i.e., 
where $\ad h$ is diagonalizable with integral eigenvalues. 
By Lemma~\ref{lem:zetacomplex}(ii) 
we may even assume that $\tau \zeta^2 = \id_{\g_\C}$, 
so that $\zeta^2 = \tau$ and therefore $\zeta(\g) = \g^c$. 

\begin{prop} Let $q_S \: S^U = \Gamma_{G}(C_U) \to \Gamma'_G(C_U) \subeq G_\C$ 
be the universal covering map of $\Gamma'_G(C_U)$, 
where $G_\C$ and $G$ are the $1$-connected 
Lie groups with Lie algebra~$\g_\C$ and $\g$, respectively.
Then $q_S(S^U_{\rm inv})$ is contained 
in the connected subgroup $\Fix(\tau_G \beta_{\pi i}) \subeq  
\Fix(\zeta_G^4)$ of $G_\C$. 
\end{prop}

\begin{prf} To apply Lemma~\ref{lem:zetacomplex}(ii), we simply have 
to observe that $q_S(G) = (G_\C)^\sigma$ is called $G$ in 
Lemma~\ref{lem:zetacomplex} and that 
\[ \beta_{\pi i}(q_S(s))
= q_S(\beta_{\pi i}(s)) = q_S(\tau_G(s)) = \tau_G(q_S(s)) 
\quad \mbox{ for } \quad s \in S^U_{\rm inv}.\] 
The subgroup $\Fix(\zeta_G^4) = (G_\C)^{\zeta_G^4}$ is connected 
by Theorem~\ref{thm:fixedpoint}. 
\end{prf}

\begin{rem} (a) One can even go one step further than the preceding
proposition by using the same trick as in the proof of 
Lemma~\ref{lem:autoholo}: Let $g \in S^U_{\rm inv} \subeq G$ 
and consider the corresponding analytic extension 
\[ \beta^g \: \oline{\cS_\pi} \to S^U, \quad z \mapsto \beta_z(g) \] 
of the orbit map of~$g$. Then the argument in the proof of 
Lemma~\ref{lem:autoholo} shows that 
$\beta^g(\oline{\cS_\pi}) \subeq g \Gamma_N(C_U)$ for 
$\fn = C_U - C_U$, so that we obtain 
in particular 
\[ \fn \ni \frac{d}{dt}\Big|_{t = 0}  g^{-1}\beta^g(t) = \Ad(g)^{-1}h - h.\] 
We conclude that 
\[ \Ad(g) h \in h  + \fn \quad \mbox{ for } \quad g \in S^U_{\rm inv}.\] 
Therefore $S^U_{\rm min}$ is contained in a Lie subgroup 
$B \subeq G$ satisfying 
\[ \Ad(B) h - h \subeq  \fn \quad \mbox{ and } \quad 
\eta_G(B) \subeq \Fix(\tau_G \beta_{\pi i}) \subeq \Fix(\zeta_G^4).\]
For the Lie algebra $\fb$ of $B$ this implies that 
$[\fb,h] \subeq \fn$, so that the semisimplicity of $\ad h$ yields 
\[ \fb \subeq \fn + \fb_0(h) \subeq \fn + \fh_0(h),\] 
where the last equality follows from the equality of 
$e^{\pi i \ad h}$ and $\tau$ on $\fb$. 

As $\fh_0 \subeq \L(S_V) = \L(S^U_{\rm inv})$ and 
the corresponding integral subgroup $H_0(h) \subeq G$ is contained 
in $S^U_{\rm inv}$, we have 
\[ S^U_{\rm inv} \cap B \subeq (S^U_{\rm inv} \cap N) H_0(h),\] 
so that the main point is to understand the subsemigroup 
\[ S^U_{\rm inv} \cap N.\] 

(b) The subgroup $\Gamma := e^{\R \ad h} \rtimes \{\1,\tau\} \subeq \Aut(\g)$ 
is abelian and $\ad h$ is diagonalizable over $\R$. 
Its Zariski closure is generated by the single element 
$\gamma := e^{\ad h} \tau$ because $\gamma^2 = e^{2 \ad h}$ 
generates a Zariski dense subgroup of $e^{\R \ad h}$. 
Hence Theorem~\ref{thm:fixedpoint} implies that the subgroup 
$(G_\C)^\Gamma$ is connected. Its Lie algebra is 
$\g_\C^\Gamma = \fh_{\C,0}(h)$ and contains $\fh_0(h)$ as a real form. 
\end{rem}

Each automorphism $\beta_z \in \Aut(G_\C)$ commutes with the 
holomorphic involution $\tau$, and hence with the holomorphic 
antiinvolution $g^\sharp = \tau(g)^{-1}$. 
As $G$ and $S^U$ are invariant under $\sharp$ because 
$C_{U}$ is invariant under $-\tau$, it follows that 
$S^U_{\rm inv}$ is $\sharp$-invariant as well. 
Therefore $g = h \exp(x) \in \zeta_G(S^U_{\rm inv}) \subeq (S^U)^{\tau_S}$ 
implies that $\exp(2x) = g^\sharp g\in \zeta_G(S_{\rm inv})$. 
But it is not clear if this implies that 
$\exp(x) \in \zeta_G(S^U_{\rm inv})$. 

The following question is 
of a similar nature. 
Let $x \in \g$ and suppose that $z := e^{yi \ad h}x \in \g_\C$ satisfies 
$\exp 2z \in \Gamma'_G(W) \subeq G_\C$. Does this imply 
that $\exp(z) \in \Gamma'_G(W)$? It seems that such questions 
are hard to answer,  as the following example shows.

\begin{ex} \mlabel{ex:5.7}
Consider the subsemigroup 
\[ S := \{ g \in \GL_n(\C) \: \|g\| \leq 1\}
= \U_n(\C) \exp(C), \quad \mbox{ where } \quad 
C = \{ X \in \Herm_n(\C) \: X \leq 0\}.\] 
We consider matrices of the form 
\[ s := \|g\|^{-1} g \quad \mbox{ for } \quad 
g = \pmat{ \eps & 1 \\ 0 & \eps}, \quad \eps > 0.\] 
Then $\|s\| \leq 1$, so that $s \in S^U$. 
Moreover, $\eps^{-1}g$ is unipotent with 
\[  X := \log g = (\log \eps)\1 + \pmat{0 & \eps^{-1} \\ 0 & 0}.\] 
Then 
\[ Y := \log s = X - (\log\|g\|) \1 
= \log(\eps\|g\|^{-1})\1 + \pmat{0 & \eps^{-1} \\ 0 & 0}\] 
satisfies $s = e^Y \in S$. That $e^{tY} \in S$ holds for all $t \geq 0$ 
is equivalent to $Y$ being dissipative, i.e., to 
\[ 0 \geq \shalf(Y + Y^*) 
= \log(\eps\|g\|^{-1})\1 + \frac{1}{2\eps}  \pmat{0 & 1\\ 1 & 0}\] 
(Remark~\ref{rem:contr-semigroup}(d)), 
which is equivalent to 
$\log(\eps) - \log(\|g\|) + \frac{1}{2\eps} \leq 0.$
For $\eps \to 0$, we have $\|g\| \to 1$, 
and $\frac{1}{2\eps} > -\log(\eps)$ if $\eps$ is sufficiently small. 
For any such $\eps$, we then have $Y \not\in \L(S)$, 
although $e^Y \in S$. 
\end{ex}

\subsection{Extensions to infinite dimensions}

\subsubsection{Wick rotations for non-uniformly continuous actions} 

It would be very interesting to understand to which extent Section~\ref{sec:3} 
can be generalized beyond uniformly continuous actions on Banach spaces, 
including $W^*$-dynamical systems. A natural setting would be that 
$E$ is a Banach space, endowed with the following data: 
\begin{itemize}
\item A continuous involution $\tau\in \GL(E)$; we write $E = E^+ \oplus E^-$, 
$E^\pm := \ker(\tau\mp \1)$ for the 
$\tau$-eigenspace decomposition.  
\item A subspace $E_* \subeq E'$ of the topological dual space 
which is {\it norm-determining} in the sense that 
$\|v\| = \sup \{ |\alpha(v)| \: \alpha \in E_*, \|\alpha\| \leq 1\}$. 
\item An $\R$-action 
$\alpha \: \R \to \GL(E)$ commuting with $\tau$ such that, 
for every $\lambda \in E_*$ and for every $v \in E$, the functions 
$t \mapsto \lambda(\alpha_t(v))$ are continuous. We say that $\alpha$ is 
$E_*$-{\it weakly continuous}. 
\item A pointed closed convex cone $W \subeq E^c := E^+ + i E^- \subeq E_\C$, 
invariant under the complex linear extension of 
$-\tau$ and the one-parameter group $(\alpha_t)_{t \in \R}$. 
\end{itemize}

We say that, for $v \in E_\C$ and $z_0 \in \C$, the element  
$\alpha^v(z_0) \in E_\C$ exists, if the orbit map $\alpha^v(t) := \alpha_t(v)$ 
extends analytically to an $E_*$-weakly continuous map on a closed strip 
$\oline{\cS_{a,b}}$ containing $z_0$. 

We expect a natural analog of Lemma~\ref{lem:4.13a} to hold. 
If $x  \in E$ is such that $\alpha_{\pi i}(x)$ exists and equals $\tau(x)$, 
then we should have an $E_*$-weakly convergent expansion 
$x = \sum_{n \in \Z} x_n$ 
with $\alpha_t(x_n) = e^{tn} x_n$ for $t \in \R$ and $\tau(x_n) = (-1)^n x_n$. 
This reduces the interesting situations to the case where 
$\zeta := \alpha_{\pi i/2}$ exists on a dense subspace 
and satisfies $\zeta^4 =\1$. 
As we cannot expect the expansion of $x$ 
to be finite, the arguments in the proof of Proposition~\ref{prop:keylema} 
fail. Presently, we are not aware of examples, where 
the conclusion of Proposition~\ref{prop:keylema} fails. 

As Olshanski semigroups and the extension of unitary representations 
also works to some extent for Banach--Lie groups 
\cite{MN12}, one may expect that large portions of our results 
can be generalized to Banach--Lie groups endowed with a 
suitably continuous action of $\R^\times$, encoding the modular objects. 

\subsubsection{The subsemigroup $S_V \subeq \U(\cH)$} 

It would be nice to find suitable regularity properties 
of $V$ that guarantee that the subsemigroup $S_V = \{ g \in \U(\cH) \:  gV \subeq V\}$ 
in the full unitary group is large in some sense. 
Of course, one could assume that it has interior points, 
but that this never leads to proper subsemigroups is easy to see: 

\begin{prop} \mlabel{prop:5.8}
Let $O \subeq \U(\cH)$ be an open subset. Then there exists an 
$N \in \N$ such that 
\[ O^N  = \{ g_1 \cdots g_N \: g_j \in O\} = \U(\cH).\] 
In particular, every subsemigroup $S \subeq \U(\cH)$ with interior points coincides
with~$\U(\cH)$.
\end{prop}

\begin{prf} Since the exponential function 
$\exp \:  \fu(\cH) = \{ X \in B(\cH) \: X^* = - X\}\to \U(\cH)$ is surjective, the open subset 
$\exp^{-1}(O)$ is non-empty. Using spectral calculus, we find an 
$n \in \N$ and an element $X \in \exp^{-1}(O)$ such that 
$\Spec(X) \subeq \frac{2\pi i}{n} \Z$. 
Then $g := e^X \in O$ is of finite order~$n$. Hence $\1 \in O^n$.

Let $B_r\subeq \fu(\cH)$ be the open operator ball of radius 
$r$ with center $0$. 
Pick $m \in \N$ such that $\exp(B_{\leq \pi/m}) \subeq O^n$.  
Then $(O^n)^m \supeq \exp(B_{\leq \pi}) = \U(\cH)$. 
\end{prf}

\appendix

\section{Conjugation with unbounded operators} 
\mlabel{app:a} 

The following proposition provides a direct path to the main 
ingredients of the Araki--Zsid\'o Theorem (Theorem~\ref{thm:lo3.18}), 
namely the implication 
(iii) $\Rarrow$ (iv). We need its corollary in the proof of 
Proposition~\ref{prop:2.8}. For the sake of completeness, we also include 
a proof of the Araki--Zsid\'o Theorem in this appendix. 

\begin{prop} \mlabel{prop:beta-ext}
Let $H = H^*$ be a selfadjoint operator and 
$U_t = e^{itH}$ denote the corresponding unitary one-parameter group. 
Fix $\beta > 0$. If $A \in B(\cH)$ is such that 
$A \cD(e^{-\beta H}) = A \cR(e^{\beta H})\subeq \cD(e^{-\beta H})$ and 
the operator $A_\beta := e^{-\beta H} A e^{\beta H}$ 
on $\cD(e^{-\beta H})$ extends to a bounded operator on $\cH$, 
then the following assertions hold: 
\begin{itemize}
\item[\rm(i)] The map $\alpha^A \: \R \to B(\cH), \alpha^A(t) 
:= U_t  A U_t^*$ extends to a bounded 
strongly continuous function on the closed strip 
$\oline{\cS_\beta} = \{ z \in \C \: 0 \leq \Im z \leq \beta\}$ which is 
holomorphic on $\cS_\beta$  and satisfies $\alpha^A(\beta i) = A_\beta$.
\item[\rm(ii)] $\|\alpha^A(z)\| \leq \max(\|A\|,\|A_\beta\|)$ 
for $z \in \oline{\cS_\beta}$
\item[\rm(iii)] $\alpha^A(z + t) = U_t \alpha^A(z) U_t^*$ 
for $z \in \oline{\cS_\beta}, t \in \R$. 
\end{itemize}
\end{prop}

\begin{prf} Let $\cH_{\rm fin} \subeq \cH$ denote the 
dense subspace of vectors contained in spectral subspaces for 
$H$ corresponding to bounded intervals. 
Let $\xi, \eta \in \cH_{\rm fin}$, so that both are entire vectors of 
exponential growth for $(U_t)_{t \in \R}$. 
Then \[ \alpha^{\xi,\eta}(z) 
:= \la e^{-i\oline z H}\xi, A e^{-i z H}\eta \ra \] 
is an entire function with 
$\alpha^{\xi,\eta}(t) = \la \xi, \alpha^A(t)\eta \ra$ and 
\begin{equation}
  \label{eq:3.17a}
\alpha^{\xi,\eta}(t + \beta i)
= \la  e^{-\beta H} U_{-t} \xi, A e^{\beta H} U_{-t} \eta \ra 
= \la  U_{-t} \xi, A_\beta U_{-t} \eta \ra 
\quad \mbox{ for } \quad t \in \R.
\end{equation}
From \cite[Thm.~12.9]{Ru86} we now derive that 
\begin{equation}
  \label{eq:formesti-a}
\|\alpha^{\xi,\eta}(z)\| \leq \max(\|A\|,\|A_\beta\|) \cdot \|\xi\| \|\eta\| \quad \mbox{ for } \quad 
z \in \oline{\cS_\pi}
\end{equation}
because this estimate holds on $\partial \cS_\beta = \R \cup (\beta i + \R)$. 
The map 
\[ \cH_{\rm fin} \times \cH_{\rm fin} \to  \cO(\cS_\pi), \quad 
(\xi,\eta) \mapsto \alpha^{\xi,\eta} \] 
is sesquilinear, and continuous with respect to the sup-norm on $\cO(\cS_\pi)$ 
by \eqref{eq:formesti-a}, hence it 
extends to a continuous map on $\cH \times \cH$ 
because $\cH_{\rm fin}$ is dense in $\cH$. 
From the one-to-one isometric correspondence between bounded 
operators and continuous sesquilinear maps on $\cH$ via 
\begin{equation}
  \label{eq:fadef}
\alpha^{\xi,\eta}(z) = \la \xi, \alpha^A(z) \eta \ra \quad \mbox{ for } \quad \xi,\eta \in\cH,
\end{equation}
we thus obtain a weakly continuous bounded map 
$\alpha^A \: \oline{\cS_\beta} \to B(\cH)$ which is weakly holomorphic 
on~$\cS_\beta$. That the function $\alpha^A \: \cS_\beta \to B(\cH)$ is holomorphic 
follows from \cite[Cor.~A.III.5]{Ne00}. 
It remains to show that it is strongly continuous on $\oline{\cS_\beta}$, 
which is done below. 

\nin (ii) follows from \eqref{eq:formesti-a}. 

\nin (iii) follows by analytic continuation because it holds for $z \in \R$. 

\nin (i) (continued) For $\eta \in \cH$, we 
consider the functions $\alpha^{A,\eta} \: \oline{\cS_\beta} \to \cH, z \mapsto \alpha^A(z)\eta$. 
By \eqref{eq:formesti-a}, we have 
$\|\alpha^{A,\eta}\|_\infty \leq \max(\|A\|,\|A_\beta\|) \|\eta\|$, so that the map 
\[ \cH \to \ell^\infty(\oline{\cS_\beta}, \cH), \quad \eta \mapsto \alpha^{A,\eta} \] 
is linear and continuous. Hence it suffices to verify the continuity of 
$\alpha^{A,\eta}$ for $\eta \in \cH_{\rm fin}$. 
For $z = x + i y \in \cS_\beta$, we have $0 \leq y \leq \beta$, so that 
\[ A e^{-izH} \eta 
\in A \cD(e^{-\beta H})
\subeq \cD(e^{-\beta H})
\subeq \cD(e^{-y H}) = \cD(e^{izH})\] 
(cf.\ \cite[Lemma~A.2.5]{NO18} for the next to last inclusion). We therefore have 
\[ \alpha^{A,\eta}(z) = e^{izH} A e^{-izH} \eta \quad \mbox{ for } \quad 
z \in \oline{\cS_\beta}.\] 
As the multiplication of operators is strongly continuous on bounded 
subsets of $B(\cH)$, (iii) shows that it suffices to verify 
the continuity of $\alpha^{A,\eta}$ on the line segment $\{ yi \: 0 \leq y \leq \beta\}$. For $0 \leq y,y_0 \leq \beta$, we have 
\begin{align} \label{eq:align}
\alpha^{A,\eta}(yi) 
&= e^{-yH} A e^{yH}\eta 
=  e^{-yH} A (e^{yH}\eta - e^{y_0H}\eta)
+ e^{-yH} A e^{y_0H}\eta.
\end{align}

Let $E$ denote the spectral measure of $H$, so that 
$H = \int_\R x \, dE(x)$. For $\xi \in \cH$ we obtain the positive 
finite measure $E^\xi := \la \xi, E(\cdot) \xi\ra$. 
Now, for $\xi \in \cD(e^{-\beta H})$, the function 
$\oline{\cS_\beta} \to \cH, z \mapsto e^{izH}\xi$ is 
continuous, because the kernel 
\[ (z,w) \mapsto 
\la e^{iwH}\xi, e^{izH}\xi \ra 
= \int_\R e^{i(z - \oline w)t}\, dE^\xi(t) \] 
is continuous on $\oline{\cS_{2\beta}}$ by 
the Dominated Convergence Theorem (\cite[Lemma~A.2.5]{NO18}). 
We conclude that the 
second summand in \eqref{eq:align} is a continuous function of~$y$.
We further have 
\[ \|e^{izH}\xi\|^2 
= \int_\R e^{- 2(\Im z) t}\, dE^\xi(t)  
\leq \max(\|\xi\|^2, \|e^{-\beta H}\xi\|^2) \] 
by the convexity of the Laplace transform of the measure $E^\xi$ 
(\cite[Prop.~V.4.3]{Ne00}). 
This implies that 
\begin{equation}
  \label{eq:doubleesti}
\|e^{- yH}\xi\| \leq \max(\|\xi\|, \|e^{-\beta H}\xi\|), 
\end{equation}
and thus 
\begin{align*}
 \|e^{-yH} A (e^{yH}\eta - e^{y_0H}\eta)\|
&\leq \|A\| \|e^{yH}\eta - e^{y_0H}\eta\| 
+ \|e^{-\beta H} A (e^{yH}\eta - e^{y_0H}\eta)\|\\
&=  \|A\| \|e^{yH}\eta - e^{y_0H}\eta\| 
+ \|A_\beta\|\|e^{-(\beta - y)H}\eta - e^{-(\beta - y_0)H}\eta\|.
\end{align*}
This estimate implies the continuity in $y_0$ 
of the first summand in \eqref{eq:align}, 
and this completes the proof of~(i).
\end{prf}

The estimate \eqref{eq:doubleesti} has an  
interesting consequence: 

\begin{cor}
\mlabel{cor:unb-cont}
Let $X$ be a topological space and $f \: X \to \cD(e^{-\beta H})$ 
be a function. 
If the two maps $f \: X \to \cH$ and 
$e^{-\beta H} \circ f \: X \to \cH$ are continuous, 
then the composition $e^{izH} \circ f \: X \to \cH$ is continuous 
for every $z \in \oline{\cS_\beta}$. 
\end{cor}

\begin{thm} \mlabel{thm:lo3.18-app} 
{\rm(Characterization of $V$-real operators)} 
For $A \in B(\cH)$, the following are equivalent: 
  \begin{itemize} 
  \item[\rm(i)] $A \in \cA_V$, i.e., $A V \subeq V$. 
  \item[\rm(ii)] $A^*V' \subeq V'$. 
  \item[\rm(iii)] $J_VA^*J_V \in \cA_V$. 
  \item[\rm(iv)] $J_VAJ_V  \Delta_V^{1/2} \subeq \Delta_V^{1/2} A$. 
  \item[\rm(v)] $\Delta_V^{1/2} A \Delta_V^{-1/2}$ is defined on 
$\cD(\Delta_V^{-1/2})$ and coincides there with $J_V A J_V$. 
  \item[\rm(vi)] The map $\alpha^A \: \R \to B(\cH), \alpha^A(t) 
= \Delta_V^{-it/2\pi} A \Delta_V^{it/2\pi}$ extends to a bounded 
strongly continuous function $\alpha^A$ on the closed strip 
$\oline{\cS_\pi} = \{ z \in \C \: 0 \leq \Im z \leq \pi\}$ which is 
holomorphic on $\cS_\pi$ and satisfies $\alpha^A(\pi i) = J_VAJ_V$. 
  \end{itemize}
If these conditions are satisfied, then 
\begin{itemize}
\item[\rm(a)] $\|\alpha^A(z)\| \leq \|A\|$ for $z \in \oline{\cS_\pi}$
\item[\rm(b)] $\alpha^A(z + t) = \Delta_V^{-it/2\pi} \alpha^A(z) \Delta_V^{it/2\pi}$ 
for $z \in \oline{\cS_\pi}, t \in \R$. 
\item[\rm(c)] $\alpha^A(\oline z + \pi i) = J_V  \alpha^A(z) J_V$ 
for $z \in \oline{\cS_\pi}$. 
\item[\rm(d)] $\alpha^A(t)V \subeq V$ and 
$\alpha^A(t + \pi i)V' \subeq V'$ for all $t \in \R$. 
\end{itemize}
\end{thm}

\begin{prf} (i) $\Rarrow$ (ii): 
If $A V \subeq V$ and $v  \in V$, $w \in V'$, then 
$\Im \la A^* w, v \ra = \Im \la w, Av \ra = 0$ shows that 
$A^*V' \subeq V'$. 

\nin (ii) $\Rarrow$ (i) follows by apply applying the implication ``(i) $\Rarrow$ (ii)'' 
to $V'$ and $A^*$ and using that $A = (A^*)^*$ and $V = (V')'$. 

\nin (ii) $\Leftrightarrow$ (iii): From $V' = J_V V$, it follows that 
$A^*V' \subeq V'$ is equivalent to $A^* J_VV \subeq J_VV$, which is (iii). 

\nin (i) $\Leftrightarrow$ (iv): For the antilinear involution $\sigma_V 
= J_V \Delta_V^{1/2}$, condition (iv) 
is equivalent to $A \sigma_V \subeq \sigma_V A$, i.e., 
to 
\[ A \cD(\sigma_V) = A (V + i V) \subeq V + i V = \cD(\sigma_V)
\quad \mbox{ and } \quad A \sigma_V = \sigma_V A \quad \mbox{ on } \quad V.\]
 This is equivalent to (i).

\nin (i) $\Leftrightarrow$ (v): Conjugating with $J_V$, we see that (v) 
is equivalent to $\sigma_V A \sigma_V^{-1} = \sigma_V A \sigma_V$ being defined on $\cD(\Delta_V^{1/2}) 
= J_V \cD(\Delta_V^{-1/2})$ and that it equals $A$ on this space. 
This in turn is equivalent to (i). 

\nin (v) $\Rarrow$ (vi) follows  from 
Proposition~\ref{prop:beta-ext} with $H = -\frac{1}{2\beta} \log(\Delta_V)$ 
and $\Delta_V^{1/2} = e^{-\beta H}$. 

\nin (vi) $\Rarrow$ (v): If (vi) is satisfied, then \eqref{eq:3.17a} 
in the proof of Proposition~\ref{prop:beta-ext} yields 
for $\xi,\eta \in \cH_{\rm fin}$ the relation 
\begin{equation}
  \label{eq:12}
\la A \Delta_V^{-1/2} \xi, \Delta_V^{1/2} \eta \ra 
= \la J_V A J_V \xi, \eta \ra.
\end{equation}
As the dense subspace $\cH_{\rm fin}$ is a core of $\Delta_V^{-1/2}$ and 
$\Delta_V^{1/2}$, the equality \eqref{eq:12} holds 
for $\xi \in \cD(\Delta_V^{-1/2})$ and $\eta \in \cD(\Delta_V^{1/2})$. 
It follows that 
\[ \Delta_V^{1/2}  A \Delta_V^{-1/2} \xi = J_V A J_V \xi
\quad \mbox{ for } \quad \xi \in \cD(\Delta_V^{-1/2}), \] 
which is (v). \\

Now we assume that the equivalent conditions (i)-(vi) are satisfied. 
From (ii) and (iii) in Proposition~\ref{prop:beta-ext}, we get~(a) and~(b). 
For $z \in \R$, we derive (c) from (vi) and~(b), and for 
general $z \in \oline{\cS_\pi}$, it follows by analytic continuation. 
Finally, (d) follows from the invariance of $V$ under 
$\Delta_V^{it}$ for $t \in \R$ and $J_V V = V'$.
\end{prf}

\section{Some facts on Lie groups} 

\begin{lem} \mlabel{lem:expinj}
Let $G$ be a finite dimensional Lie group with Lie 
algebra $\g$ and $x,y \in \g$ with $\exp x = \exp y$. 
If $\exp$ is not singular in $x$, then $[x,y] =0$ and 
$\exp(x-y) = e$. 

If, in addition, $G$ is simply connected and $\ad x$ and $\ad y$ have 
real spectrum, then $x = y$. 
\end{lem}

\begin{prf} The first assertion follows from \cite[V.6.7]{HHL89}. 
If $\ad x$ and $\ad y$ have real spectrum, then 
$\exp$ is regular in $x$, so that $[x,y] = 0$ and $z := x-y$ satisfies 
$\exp(z) = e$. 
The latter condition implies $e^{\ad z} = \1$, so that 
$\ad z$ is semisimple with purely imaginary spectrum. 
On the other hand, $[\ad x, \ad y] = \ad [x,y] = 0$ implies that 
$\Spec(\ad z) \subeq \Spec(\ad x) - \Spec(\ad y) \subeq \R$ 
(there exists a common generalized eigenspace decomposition). 
Combining both facts, we see that $\ad z= 0$, i.e., 
$z \in \fz(\g)$. If $G$ is simply connected, then 
$\exp\res_{\fz(\g)}$ is injective because 
$Z(G)_0 = \exp(\fz(\g))$ is simply connected 
(\cite[Thm.~11.1.21]{HN12}). This implies $z = 0$. 
\end{prf}

\begin{thm} \mlabel{thm:fixedpoint}  
Let $G$ be a $1$-connected Lie group and 
let $\Gamma \subeq \Aut(G)$ be a subgroup such that the Lie algebra 
$\g$ is a semisimple $\Gamma$-module. 
Then the following assertions hold: 
\begin{itemize}
\item[\rm(i)] There exists a $\Gamma$-invariant Levi decomposition  
$G \cong R \rtimes S$, so that the subgroup of $\Gamma$-fixed points is 
$G^\Gamma \cong R^\Gamma \rtimes S^\Gamma$. 
\item[\rm(ii)] The group $R^\Gamma$ is connected. 
\item[\rm(iii)] If the action of $\Gamma$ on the Lie algebra $\fs$ of $S$ 
has a relatively compact image in $\Aut(\fs)\cong \Aut(S)$ which 
contains a dense cyclic subgroup, then $S^\Gamma$ is connected.
\begin{footnote}
{For any element $\gamma \in \Gamma$ for which $\gamma^\Z$ is dense in 
$\Gamma$ we then have the same group of fixed points. Note also that 
this assumption is satisfied if $\Gamma$ is a product of a torus 
and a finite cyclic group.}  
\end{footnote}
\item[\rm(iv)] If $\eta_S \: S \to S_\C$ is the universal complexification, 
then the $\Gamma$-action on $S$ induces an action on $S_\C$. 
If the image of $\Gamma$ in the algebraic group $\Aut(\fs)$ is 
generated by a single semisimple automorphism in the Zariski topology, 
then $(S_\C)^\Gamma$ is connected.\begin{footnote}
{In Borel's book \cite{Bor91} one finds in particular 
that centralizers of complex tori are connected 
(\cite[Cor.~11.12]{Bor91}). Since every torus 
contains a single element with the same centralizer 
(\cite[Prop.~8.18]{Bor91}) this follows from the present statement of (iv).  
}\end{footnote}

Further $\eta_S(S^\Gamma)$ is an open subgroup in the group 
$\eta_S(S)^\Gamma = (S_\C)^{\Gamma,\sigma}$, where 
$\sigma$ is the complex conjugation on $S_\C$ with fixed point set 
$\eta(S) = (S_\C)^\sigma$.  
\end{itemize}
\end{thm}

\begin{prf} (i) With \cite[Prop.~I.2]{KN96} we find a $\Gamma$-invariant 
Levi decomposition $\g = \fr \rtimes \fs$, so that we obtain 
a Levi decomposition $G \cong R \rtimes S$, where 
$R$ is solvable, $S$ is semisimple and both are $1$-connected and 
$\Gamma$-invariant. This proves (i). 

\nin (ii) We argue by induction on 
the dimension of $R$. If $R$ is abelian, 
then this $1$-connected group is isomorphic to some $\R^n$ and 
$\Gamma$ acts by  linear maps. This implies that $R^\Gamma$ is a linear subspace, 
hence connected. 

If $R$ is not abelian, then its commutator subgroup 
$R'  = (R,R)$ has smaller dimension and 
its Lie algebra $\fr' = [\fr,\fr]$ is a proper $\Gamma$-invariant 
ideal of $\fr$. Let $\fn \supeq \fr'$ be a maximal 
proper $\Gamma$-invariant ideal of $\fr$ and let $N \trile R$ be the 
corresponding normal integral subgroup. 
Since $R$ is $1$-connected, $N$ is closed and $1$-connected and   
the abelian quotient group $Q := R/N$ is also $1$-connected 
(\cite[Thm.~11.1.21]{HN12}). 
As $N$ is $1$-connected, our 
induction hypothesis implies that $N^\Gamma$ is connected. 
As $N$ is $\Gamma$-invariant, $Q$ inherits a natural $\Gamma$-action 
and since $Q$ is abelian, the above argument shows that the fixed 
point group $Q^\Gamma$ is connected. 

Clearly, $q(R^\Gamma) \subeq Q^{\Gamma}$, and we claim that we actually 
have equality. Two cases may occur. 
If $Q^{\Gamma} = \{e\}$, then $R^\Gamma = N^\Gamma$ is connected. 
If $Q^{\Gamma} \not=\{e\}$, then it is a connected 
subgroup of positive dimension. 
As the action of $\Gamma$ on $\fr$ is semisimple, there exists a 
$\Gamma$-invariant linear subspace $\fe \subeq \fr$ complementing~$\fn$. 
Then $\L(q) \: \fe \to \fq$ is a linear $\Gamma$-equivariant 
isomorphism, and since $\exp_Q \: (\fq,+) \to Q$ also is an isomorphism 
of Lie groups, it follows that 
\[ \exp_Q \circ \L(q) = q \circ \exp_R \: \fe^\Gamma \to Q^\Gamma \] 
is a bijection. Although $\fe$ may not be a Lie subalgebra 
of $\fr$, the preceding argument shows that 
$R^\Gamma/N^\Gamma \cong q(R^\Gamma) = Q^\Gamma$. 
As $N^\Gamma$ and $Q^\Gamma$ are connected, we conclude that the 
group $R^\Gamma$ is connected as well. 

\nin (iii) 
Replacing $\Gamma$, considered as a subgroup of $\Aut(\fs) \cong \Aut(S)$, 
by its compact closure does not change the subgroup of fixed points 
because the action of $\Aut(\fs) \cong \Aut(S)$ on $S$ 
is smooth (\cite[Thm.~11.3.5]{HN12}). 
So we may w.l.o.g.\ assume that $\Gamma$ is compact. 
It therefore is contained in a maximal compact subgroup $C \subeq \Aut(\fs)$ 
because $\Aut(\fs)$ is an algebraic group, hence has only finitely 
many connected components (\cite[\S 12.4]{HN12}). 

Now $C \cap \Aut(\fs)_0 = C \cap \Ad(S)$ is maximal compact in the identity 
component, and therefore $K := \{ s \in S \: \Ad(s) \in C \}$ 
is maximal compactly embedded in $S$. We conclude that 
$K$ is $1$-connected and therefore that 
$K \cong Z(K) \times (K,K)$, where $Z(K)$ is a vector space 
and $(K,K)$ is $1$-connected compact, a maximal compact subgroup of~$S$ 
(\cite[Thm.~12.1.18]{HN12}). 
As $K$ is invariant under the action of $C$ on $S$, it is in particular 
invariant under $\Gamma$. Since $\Gamma$ acts by automorphisms on 
$K$, it preserves its center $Z(K)$ and its commutator subgroup~$(K,K)$. 
Let $\fp \subeq \fs$ be the orthogonal complement 
of the Lie algebra $\fk$ of $K$ with respect to the Killing form. 
Then the polar map $K \times \fp \to S, (k,x) \mapsto k \exp x$ 
is a $\Gamma$-equivariant diffeomorphism. 
We thus obtain 
\[ S^\Gamma = K^\Gamma \exp(\fp^{\Gamma})
\cong (K')^\Gamma \times Z(K)^\Gamma \times \fp^\Gamma.\] 
As $Z(K)$ is a vector space, the group $Z(K)^\Gamma$ is a linear 
subspace, hence connected. The same is true for $\fp^\Gamma$. 
To verify the connectedness of $(K')^\Gamma$, 
we recall that there exists a single element $\gamma \in \Gamma$ 
for which the cyclic subgroup $\gamma^\Z$ is dense in $\Gamma$, 
considered as a subgroup of $\Aut(\fs)$. 
As $\Aut(\fs) \cong \Aut(S)$ acts smoothly on $S$ 
(\cite[Thm.~11.3.5]{HN12}), 
it follows that $\Gamma$ and $\gamma$ have the same fixed points. 
Now the $1$-connectedness of the compact group $K'$ 
implies that $(K')^\Gamma = (K')^\gamma$ is connected 
(\cite[Thm.~12.4.26]{HN12}). This shows that $S^\Gamma$ is connected. 

(iv) Let $\gamma \in \Gamma\subeq \Aut(\fs)$ be a semisimple element 
for which $\Gamma$ is contained in the Zariski closure of the cyclic 
subgroup $\gamma^\Z$. Since the action of the algebraic group 
$\Aut(\fs)$ on the algebraic group $S_\C$ is algebraic, 
$\gamma$ and $\Gamma$ have the same fixed point group. 
As the group $S_\C$ is  $1$-connected, the 
connectedness of $S_\C^\gamma = S_\C^\Gamma$ 
now follows from \cite[Thm.~4.4.9, p.~214]{OV90}. 
The remaining assertions are clear.
\end{prf}

From Theorem~\ref{thm:fixedpoint}(i)-(iii), we obtain in particular: 

\begin{cor} \mlabel{cor:fixedpoint}  
Let $G$ be a $1$-connected Lie group and 
$\phi \in \Aut(G)$ an automorphism of finite order. 
Then the subgroup $G^\phi = \{ g \in G \: \phi(g) = g\}$ 
of fixed points is connected. 
\end{cor}

\section*{Acknowledgments}
We thank Yoh Tanimoto and Roberto Longo for an invitation 
to a research visit in Rome and for many discussions 
with them, Vincenzo Morinelli and Yoshimichi Ueda 
on standard subspaces and modular theory of von Neumann algebras. 
In particular, we thank Yoh Tanimoto for pointing out an inaccuracy in an earlier 
version of this paper. 

Last, but not least, we also thank Daniel Oeh and Jan Frahm for reading 
earlier versions of this manuscript.

\end{document}